\documentclass{article}

\usepackage{graphicx}
\usepackage{amssymb}
\usepackage{color}
\usepackage{wrapfig}
\usepackage{subfig}
\usepackage{url}

\usepackage[margin=1.0in]{geometry}

\newtheorem{theorem}{Theorem}

\newtheorem{prop}[theorem]{Proposition}

\newtheorem{corollary}[theorem]{Corollary}

\def\squareforqed{\hbox{\rlap{$\sqcap$}$\sqcup$}}
\def\qed{\ifmmode\else\unskip\quad\fi\squareforqed}
\def\smartqed{\def\qed{\ifmmode\squareforqed\else{\unskip\nobreak\hfil
\penalty50\hskip1em\null\nobreak\hfil\squareforqed
\parfillskip=0pt\finalhyphendemerits=0\endgraf}\fi}}

\begin{document}
	
\title{Geometric auxetics}
\author{\setcounter{footnote}{0}%
\def\thefootnote{\arabic{footnote}}
Ciprian S. Borcea$^1$ and
Ileana~Streinu$^2$ 
}
	  
\date{}

\maketitle


{
\setcounter{footnote}{0}
\def\thefootnote{\arabic{footnote}}
\footnotetext[1]{ Department of Mathematics, Rider University, Lawrenceville, NJ 08648, USA.
\url{borcea@rider.edu}
}
\footnotetext[2]{Department of Computer Science, \  Smith College,
Northampton, \ MA 01063, USA. \  \url{istreinu@smith.edu}. 
}
}

\begin{abstract}
We formulate a mathematical theory of auxetic behavior based on one-parameter deformations of periodic frameworks. Our approach is purely geometric, relies on the evolution of the periodicity lattice and works in any dimension.  We demonstrate its usefulness by predicting or \ recognizing, without \ experiment, computer simulations or \ numerical \ approximations, the auxetic capabilities of several well-known structures available in the literature. We propose new principles of auxetic design and rely on the stronger notion of expansive behavior to provide an infinite supply of planar auxetic mechanisms and several new three-dimensional structures. 
\end{abstract}


\section*{Introduction}
 
The notion of {\em auxetic behavior} emerged from renewed interest in materials with negative Poisson's ratios \cite{ENHR,GGLR,L}. Auxeticity is often introduced in suggestive terms as a property of materials which undergo a lateral widening upon stretching.
In theoretical terms, Poisson's ratio is defined in elasticity theory and involves physical properties of the material under consideration \cite{G}. Quoting from \cite{TC}: ``Poisson's ratio $\nu(\mathbf{n},\mathbf{m})$ of an elastic solid for any two specified orthogonal unit vectors $\mathbf{n}$ and $\mathbf{m}$ is the ratio of the lateral contraction in the direction $\mathbf{m}$ to the axial extension in the direction $\mathbf{n}$ due to a uniaxial tension of the material along the direction $\mathbf{n}$."

The main purpose of this paper is to show that, for crystalline materials and  man-made structures modeled as periodic bar-and-joint frameworks, a {\em purely geometric approach} to auxetic properties can be defined and studied mathematically, with no admixture of physical assumptions. 

Geometrical underpinnings for auxetic behavior have been suggested in several contexts. In dimension two, examples have \ been \ based \ on so-called {\em reentrant honeycombs} \cite{Kol,PL}, {\em missing rib models} \cite{GRSGE}, {\em rotating rigid units} \cite{GAE,G-E,SSML} and certain regular {\em plane tessellations} \cite{GS,M-ST,MRMST,M-G}.
In dimension three, auxetic behavior was proposed for $\alpha$-cristobalite and $\alpha$-quartz \cite{AE,KC,PG,YHWP}. Pairs of orthogonal directions $(\mathbf{n},\mathbf{m})$ with a negative Poisson's ratio are fairly common in cubic crystals \cite{Bau,BauSZS,N}. Considerations related to composite materials appear in \cite{M1}.
See also \cite{EL,G-G,Koe,M2}. In spite of this diversity of examples, no geometric principle underlying the auxetic properties of periodic structures  has been proposed so far.
 
\medskip
In  this paper we define a general notion of {\em auxetic path} or {\em auxetic one-parameter deformation} for periodic bar-and-joint frameworks in Euclidean spaces of arbitrary dimension $d$. This definition takes into account two fundamental aspects: (i) when deformed {\em as a periodic  structure,} the framework carries a natural reference system in its lattice of periods and (ii) in a one-parameter deformation, two instances at times $\tau_1 < \tau_2$, can be compared via the linear transformation which takes the period lattice of the latter to the period lattice of the former (consistent with continuity). We take as {\em the} essential feature of an auxetic deformation, the property of this linear transformation to be a contraction for any pair of time parameters $\tau_1 < \tau_2$.

\medskip  
For a rigorous  mathematical treatment, we rely on the deformation theory of periodic frameworks introduced in \cite{BS1} and developed in \cite{BS2,BS3}. In Section~\ref{sec:deformations} we review the necessary definitions. The point of view presented in \cite{BS3} 
implicates directly and naturally the {\em positive semidefinite cone} in the space of $d\times d$ symmetric matrices \cite{B,Gr}. By recording the Gram matrix of an independent set of generators for the periodicity lattice, a deformation path determines a parametrized curve inside this cone. We prove in Section~\ref{sec:criterion}, Theorem~\ref{thm:tangent}, that the auxetic character of a deformation path is equivalent to the fact that all tangent vectors along this  curve of Gram matrices belong (as free vectors) to the positive semidefinite cone. This equivalent version offers a good analogy between auxetic trajectories and causal lines in special relativity \cite{Borel,LEMW}. Actually, for $d=2$, this analogy is a literal matching, since the three-dimensional space of $2\times 2$ symmetric matrices carries a natural Minkowski metric. We elaborate on this topic in Section~\ref{sec:Minkowski}.

\medskip
We demonstrate the usefulness and simplicity of our proposed criterion for auxeticity by revisiting, in Section~\ref{sec:silica}, the classical tilting scenarios for the $\alpha$ to $\beta$ phase transitions in the silica polymorphs quartz and cristobalite. These one-parameter deformations are auxetic in our geometric sense.

\medskip
In Section~\ref{sec:planar}, we focus on the planar case and engage at the same time the {\em design} problem for frameworks with auxetic capabilities. We emphasize in this context the far-reaching role of the stronger notion of {\em expansive behavior} and its correlate, the concept of periodic pseudo-triangulation.  A one-parameter deformation is expansive when {\em all} distances between pairs of vertices (i.e. joints) increase or stay the same. The understanding of expansive deformations in dimension two is not an elementary matter. It relies on results initially discovered  for finite linkages \cite{S1,S2} and a recently established periodic version of a classical theorem of Maxwell \cite{BS4,Max}. Periodic pseudo-triangulations explain in full the phenomenon of two-dimensional expansive behavior.  Given that {expansive implies auxetic}, we obtain, via periodic pseudo-triangulations, an infinite series of auxetic planar designs. 
 
\medskip
In Section~\ref{sec:scenario} we return to the general case of dimension $d$ and show that, by combining our deformation theory of periodic frameworks \cite{BS1,BS3} with the auxetic principle presented in this paper, we obtain precise information about all possible auxetic trajectories. At the infinitesimal level, there is a well defined and  accessible {\em cone of infinitesimal auxetic deformations}. From a computational perspective, we find here a direct link with {\em semidefinite programming} and recent research on {\em spectrahedra} \cite{ORSV,PK,Vin}.  Then, we illustrate the main scenario with a three dimensional periodic framework with four degrees of freedom and  various capabilities for expansive or auxetic deformations. The position of the polyhedral expansive cone is explicitly determined in the larger spectrahedral auxetic cone. Despite limited general results about expansive behavior in three or higher dimensions,  we present several new three-dimensional auxetic designs suggested by necessary conditions for expansive deformations. The final section collects our conclusions.

\section{Periodic frameworks and their deformations}
\label{sec:deformations}

For auxetic behavior, a periodic framework must first of all be flexible. Although structural flexibility in framework crystalline materials was observed at an early date \cite{Pauling}, interest for a self-sufficient mathematical treatment, independent of solid state physics, developed quite recently. For our present purpose, the appropriate mathematical instrument is the deformation theory of periodic frameworks presented in \cite{BS1}, with elaborations contained in \cite{BS2,BS3,BS4,BS6}. We review the main concepts in this section.

\begin{figure}[h]
\centering
 {\includegraphics[width=0.44\textwidth]{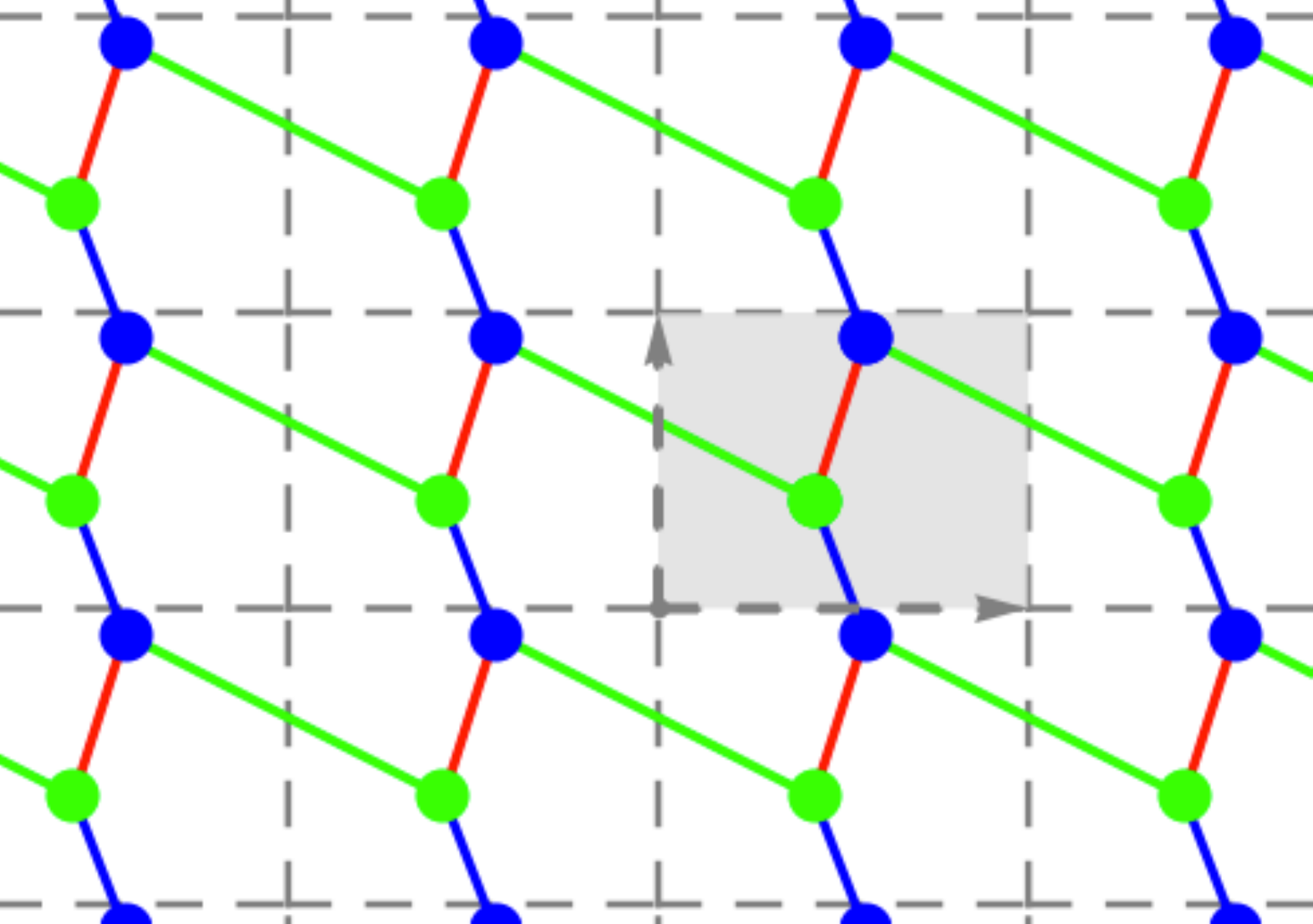}}
 \hspace{20pt}
 {\includegraphics[width=0.44\textwidth]{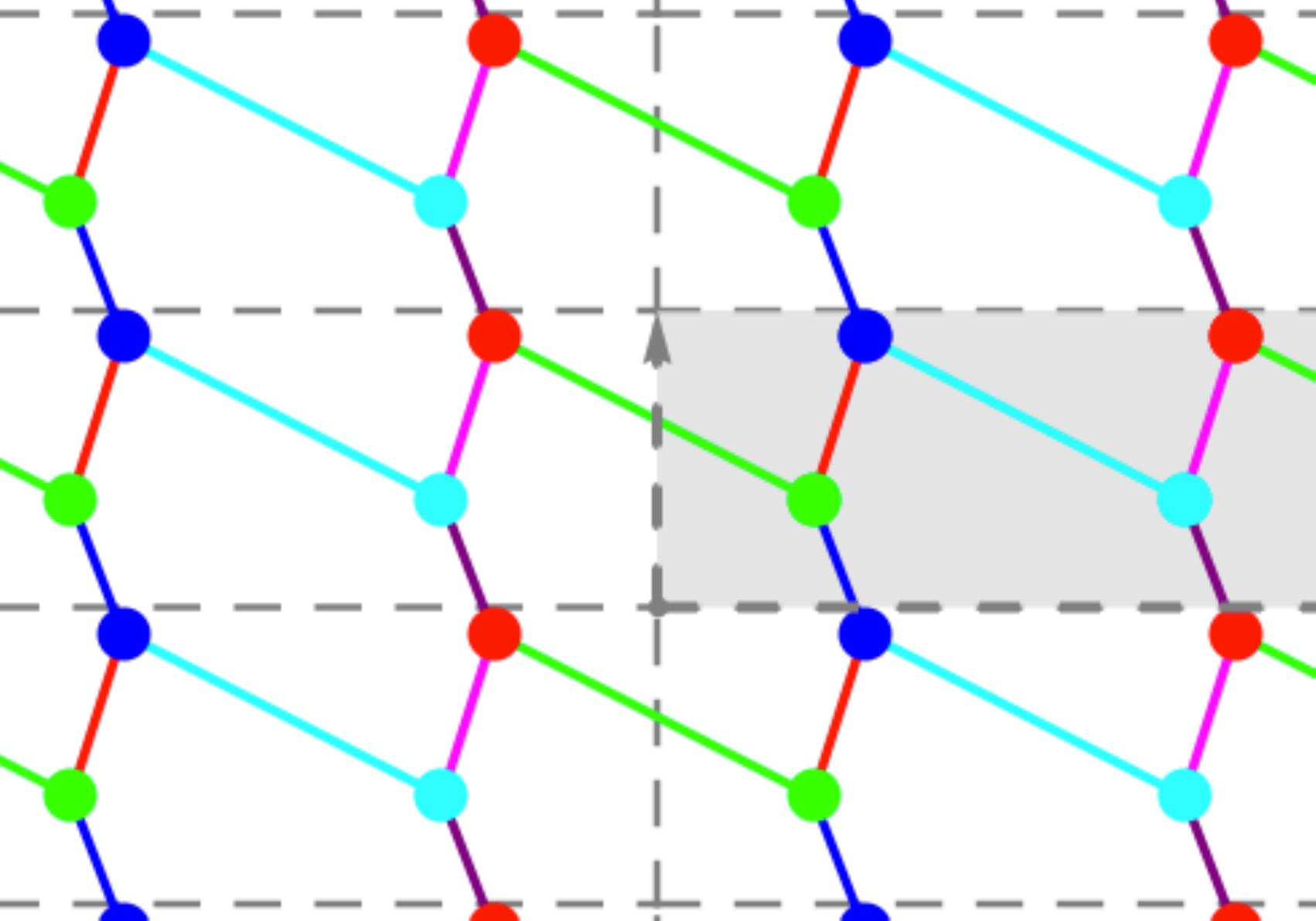}}\\
 {\includegraphics[width=0.12\textwidth]{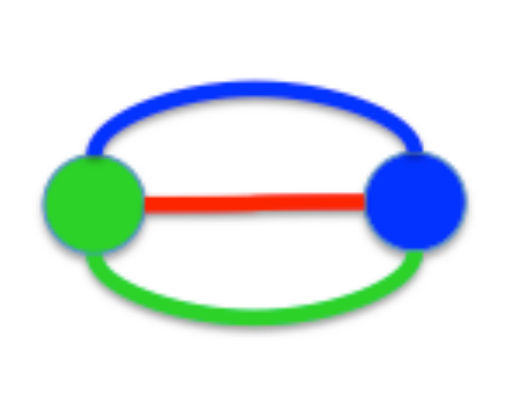}}
 \hspace{4pt}
 {\includegraphics[width=0.12\textwidth]{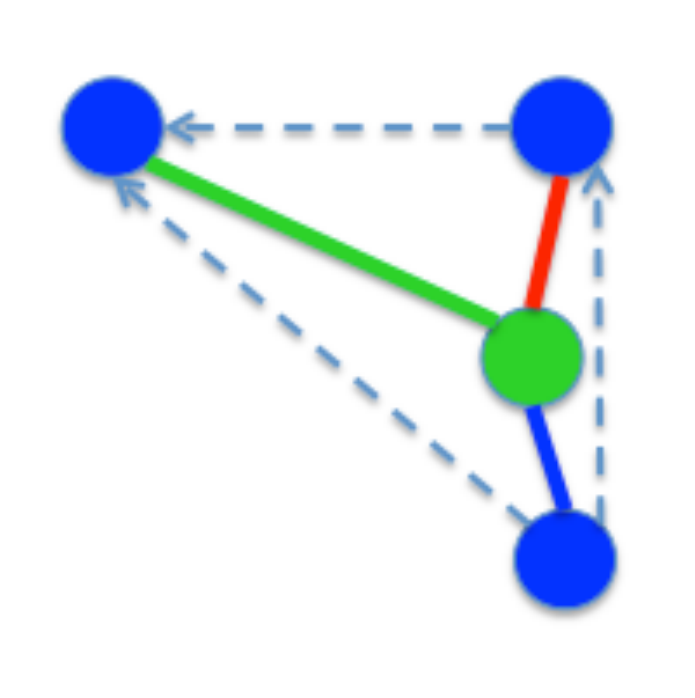}}
 \hspace{40pt}
 {\includegraphics[width=0.12\textwidth]{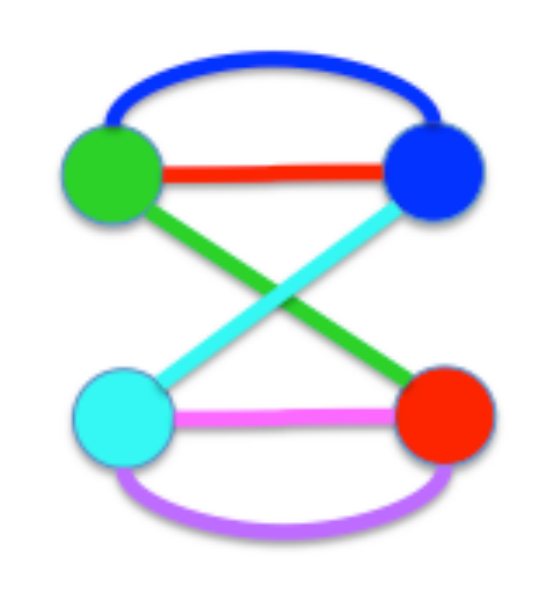}}
 \hspace{4pt}
 {\includegraphics[width=0.20\textwidth]{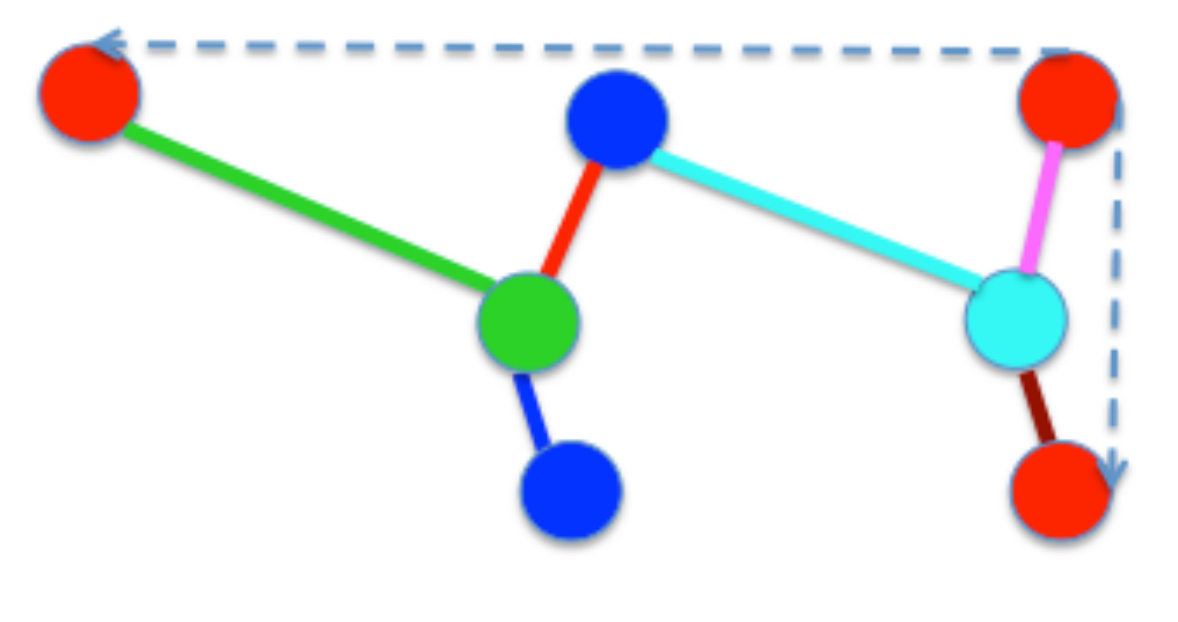}}
 \caption{ (Top) Two  $2$-periodic graphs defined on the same infinite graph $G$ and differing by the periodicity group $\Gamma$ acting on $G$. Generators of the two lattices and the corresponding unit cells are highlighted, and orbits of vertices and edges are similarly colored. (Left) For maximal periodicity, there are two vertex orbits and three edge orbits. (Right) With an index two sublattice, there are four vertex orbits and six edge orbits. (Bottom) The corresponding quotient graphs, with matching vertex and edge colors, together with a schematic representation for the placement of vertex and edge representatives.} 
 \label{FigPeriodic}
 \vspace{-14pt}
\end{figure}

\noindent
{\bf Periodic graph.}
A $d$-periodic graph is a pair $(G,\Gamma)$, where $G=(V,E)$ is a simple infinite graph with vertices $V$, edges $E$ and finite degree at every vertex, and $\Gamma \subset Aut(G)$ is a free Abelian group of automorphisms which has rank $d$, acts without fixed points and has a finite number of vertex (and hence, also edge) orbits. The group  $\Gamma$ is thus isomorphic to $Z^d$ and is called the {\em periodicity group}  of the periodic graph $G$. Its elements $\gamma \in \Gamma \simeq Z^d$ are refered to as {\em periods} of $G$. We emphasize the key role played in our definition by the periodicity group $\Gamma$ and illustrate it in Fig.~\ref{FigPeriodic}: different groups acting on the same infinite graph yield distinct periodic graphs. 

\medskip
\noindent
{\bf Quotient graph.}  To any $d$-periodic graph $(G,\Gamma)$ we associate a {\em quotient graph} $G/\Gamma=(V/\Gamma,E/\Gamma)$, whose vertices and edges correspond to vertex, resp. edge orbits in $(G,\Gamma)$.  Examples are shown in Fig.~\ref{FigPeriodic}.

\medskip
\noindent
{\bf Periodic placement of a periodic graph.}
A periodic placement (shortly, a placement) of a $d$-periodic graph $(G,\Gamma)$ in $R^d$ is defined by two functions:
\[ p:V\rightarrow R^d \ \ \mbox{and} \ \ \pi: \Gamma \hookrightarrow {\cal T}(R^d) \]

\noindent 
where $p$ assigns points in $R^d$ to the vertices $V$ of $G$ and $\pi$ is a faithful representation of the periodicity group $\Gamma$, that is, an injective homomorphism of $\Gamma$ into the group ${\cal T}(R^d)$ of translations in the Euclidean space $R^d$, with $\pi(\Gamma)$ being a lattice of rank $d$. These two functions must satisfy the natural compatibility condition:
$$ p(\gamma v)=\pi(\gamma)(p(v)) $$
 
\noindent
After choosing an independent set of $d$ generators for the periodicity lattice $\Gamma$, the image $\pi(\Gamma)$ is completely described via the $d\times d$ matrix $\Lambda$ with column vectors $(\lambda_i)_{i=1,\cdots,d}$ given by the images of the generators under $\pi$. The parallelotope spanned by these vectors defines a {\em unit cell} of the periodic placement. Examples in dimension two are shown in Fig.~\ref{FigPeriodic}.

\medskip
\noindent
{\bf Periodic framework.}
A placement which does not allow the end-points of any edge to have the same image defines a {\em $d$-periodic bar-and-joint framework} in $R^d$, with edges $(u,v)\in E$ corresponding to bars (segments of fixed length) $[p(u),p(v)]$ and vertices corresponding to (spherical) joints. Two frameworks are considered equivalent when one is obtained from the other by a Euclidean isometry. 

\medskip
\noindent
{\bf Periodic deformation.} A {\em one-parameter deformation of the periodic framework} $(G,\Gamma, p,\pi)$ is a  (smooth) family of placements  $p_{\tau}: V\rightarrow R^d$  parametrized by time $\tau \in (-\epsilon, \epsilon)$ in a small neighborhood of the initial placement $p_0=p$, which satisfies two conditions: (a) it maintains the lengths of all the edges $e\in E$, and (b) it maintains periodicity under $\Gamma$, via faithful representations $\pi_{\tau}:\Gamma \rightarrow {\cal T}(R^d)$ which {\em may change with $\tau$ and give a concomitant variation of the periodicity lattice} $\Lambda_{\tau}=\pi_{\tau}(\Gamma)$.  {Fig.~\ref{FigAnimPPT} shows a few snapshots from a one-parameter deformation of a periodic framework, and illustrates the fact that, while the abstract periodicity group $\Gamma$ continues to act on the deformed framework, its geometric counterpart $\pi(\Gamma)$ varies with time.

\medskip
\noindent
{\bf Deformation space.} Given a $d$-periodic framework $(G,\Gamma, p,\pi)$, the collection of all periodic placements in $R^d$ which maintain the lengths of all edges is called the {\em realization space} of the framework. After factoring out equivalence under Euclidean isometries, one obtains the {\em configuration space} of the framework (with the quotient topology). The {\em deformation space} is the connected component of the configuration space which contains the initial framework.

For more background regarding these concepts, the reader should consult the original papers. It is important to retain the fact that framework deformation spaces are semi-algebraic sets and notions of algebraic or differential geometry such as singularity, tangent space, dimension, apply accordingly.

\section{Auxetic one-parameter deformations}
\label{sec:criterion}

This section introduces the main definition of our paper and offers an alternative characterization. In preparation, we recall the following classical concept from operator theory.

\medskip
\noindent
{\bf Contraction operator.}
Let $T: R^d \rightarrow R^d$ be a linear operator. The {\em operator norm} (shortly, the norm) of $T$ is defined as:
$$ ||T||=sup_{|x|\leq 1} |Tx|=sup_{|x|=1} |Tx| $$

\noindent
$T$ is called a {\em contraction operator} (shortly, a contraction), when $||T||\leq 1$, and a {\em strict contraction} if $||T|| < 1$.
From $|Tx|\leq ||T||\cdot |x|$ it follows that contraction operators are characterized by the property of taking the unit ball to a subset of itself. 

\medskip
Given a one-parameter deformation $(G,\Gamma, p_{\tau},\pi_{\tau})$, $\tau\in (-\epsilon,\epsilon)$ of a periodic framework in $R^d$, the corresponding {\em one-parameter} family of periodicity lattices $\Lambda_{\tau}=\pi_{\tau}(\Gamma)$ offers a way to compare any two sequential moments $\tau_1 < \tau_2$ by looking at the {\em unique linear operator} $T_{\tau_2\tau_1}$ which takes the lattice at time $\tau_1$ to the lattice at time $\tau_2$, i.e. is defined by:
\begin{equation}\label{eq:t1t2}
 \pi_{\tau_1}=T_{\tau_2\tau_1}\circ \pi_{\tau_2} 
\end{equation}
%

\noindent
{\bf Main definition: auxetic path.} {\em A differentiable one-parameter deformation
$$(G,\Gamma, p_{\tau},\pi_{\tau})_{\tau\in (-\epsilon,\epsilon)}$$
of a periodic framework in $R^d$ is said to be an {\em auxetic path} (shortly, auxetic), when for any $\tau_1 < \tau_2$, the linear operator $T_{\tau_2\tau_1}$ defined by (\ref{eq:t1t2}) is a contraction.}

\medskip
\noindent
{\bf Remark.} This formulation provides the most intuitive connection between our geometric approach to auxetic behavior and the conventional approach based on negative Poisson's ratios. Indeed, in a contraction {\em any} vector is mapped to a vector of smaller or equal length; hence shrinking in {\em any} particular direction entails `lateral' shrinking as well, as in the formulation via Poisson's ratios.

\medskip \noindent
It is important to obtain an infinitesimal version of this criterion, that is, a characterization in terms of infinitesimal deformations along the path. This can be derived as follows. After choosing an independent set of generators for the periodicity lattice $\Gamma$, the image $\pi_{\tau}(\Gamma)$ is completely described via the $d\times d$ matrix $\Lambda_{\tau}$ with column vectors given by the images of the generators under $\pi_{\tau}$. The associated {\em Gram matrix} is given by:
$$ \omega_{\tau}=\omega(\tau)=\Lambda^t_{\tau}\Lambda_{\tau}. $$

\begin{theorem}\label{thm:tangent}
A deformation path $(G,\Gamma, p_{\tau},\pi_{\tau}), \tau\in (-\epsilon,\epsilon)$ is auxetic if and only if the curve of Gram matrices $\omega(\tau)$ defined above has all its
tangents in the cone of positive semidefinite symmetric $d\times d$ matrices.
\end{theorem}

\noindent
{\em Proof:}\
In one direction, we assume that $(G,\Gamma, p_{\tau},\pi_{\tau}), \tau\in (-\epsilon,\epsilon)$ is an auxetic path. We have to prove that:
\begin{equation}\label{eq:at0}
\dot{\omega}_0=\frac{d\omega}{d\tau}(0)
\end{equation}
\noindent
is positive semidefinite. To simplify the notation we rewrite (\ref{eq:t1t2}) as:
$$ T_{\tau,0}=T_{\tau}=\Lambda_0\Lambda_{\tau}^{-1} \ \ \mbox{or} \ \  
\Lambda_{\tau}=T_{\tau}^{-1}\Lambda_0 $$
\noindent 
which gives:
\begin{equation}\label{eq:Gram}
\omega_{\tau}=\Lambda_0^t(T_{\tau}^{-1})^tT_{\tau}^{-1}\Lambda_0.
\end{equation}
Further simplifying the notation with $D_{\tau}=T_{\tau}^{-1}$ and $\dot{D}_0=\frac{dD}{d\tau}(0)$, we see that it suffices to show that:
\begin{equation}\label{eq:enough}
\langle \dot{D}_0 x, x\rangle +\langle x, \dot{D}_0x\rangle \geq 0
\end{equation}
\noindent
This will follow from the fact that the function of $\tau$ given by 
$\langle D_{\tau}x,D_{\tau}x\rangle$ is non-decreasing. 
Indeed, since $T_{\tau_2\tau_1}$ is a contraction for $\tau_2 > \tau_1$, we have:
\begin{equation}\label{eq:increasing}
\langle D_{\tau_1}x,D_{\tau_1}x\rangle = |D_{\tau_1}x|^2=|T_{\tau_2\tau_1} D_{\tau_2}x|^2
\leq |D_{\tau_2}x|^2
\end{equation}
\noindent
Taking the derivative at $\tau=0$ in ({\ref{eq:enough}), we obtain the desired conclusion:
\begin{equation}\label{eq:derivation}
\frac{d}{d\tau} |D_{\tau}x|^2|_{\tau=0}=\frac{d}{d\tau} \langle x,D_{\tau}^tD_{\tau}x\rangle |_{\tau=0}=\langle x, (\dot{D}_0^t+\dot{D}_0)x\rangle \geq 0
\end{equation}

\medskip
\noindent
For the reverse implication, let us assume that the deformation $(G,\Gamma, p_{\tau},\pi_{\tau})_{\tau}$ gives a curve of Gram matrices which has all its velocity vectors in the positive semidefinite cone. Integrating from $0$ to $\tau$, we find:
\begin{equation}\label{eq:integrate}
\omega_{\tau}-\omega_0=P \succeq 0
\end{equation}
\noindent
with notation indicating that $P$ is a positive semidefinite operator. We may assume $\tau=1$, with $\omega_1=\Lambda_1^t\Lambda_1$ and $\omega_0=\Lambda_0^t\Lambda_0$. Our task amounts to proving that $T_1=\Lambda_0\Lambda_1^{-1}$ is a contraction. 

\medskip
In fact, noting that the argument is not affected by left multiplication of $\Lambda_i$, $i=0,1$ with an orthogonal transformation, 
we may assume $\Lambda_0=R\omega_0^{1/2}$ and $\Lambda_1=\omega_1^{1/2}$, where $R$ is orthogonal and $\omega_i^{1/2}$ stands for the unique
positive square root of $\omega_i$.  We obtain:
$$
\langle T_1x,T_1x\rangle=
\langle \Lambda_0\Lambda_1^{-1}x, \Lambda_0\Lambda_1^{-1}x \rangle = $$
$$
= \langle R\omega_0^{1/2}\omega_1^{-1/2}x, R\omega_0^{1/2}\omega_1^{-1/2}x \rangle =
\langle \omega_0^{1/2}\omega_1^{-1/2}x, \omega_0^{1/2}\omega_1^{-1/2}x \rangle = $$
$$
= \langle \omega_1^{-1/2}\omega_0 \omega_1^{-1/2}x, x \rangle =
\langle \omega_1^{-1/2}(\omega_1 - P)\omega_1^{-1/2}x, x \rangle = $$
\begin{equation}\label{eq:end}
= \langle x,x \rangle -\langle P\omega_1^{-1/2}x, \omega_1^{-1/2}x \rangle \leq 
\langle x,x \rangle
\end{equation}

\noindent
Hence $T_1$ is a contraction, and this concludes the proof. \qed

\medskip 
\noindent
{\bf Commentary and interpretation.} The {\em auxeticity} criterion formulated in terms of contraction operators is probably more intuitive:  the image of the unit ball going to a subset of itself does convey a sense of coordinated shrinking (that is, coordinated {\em growth} when increasing the time parameter). On the other hand, the equivalent characterization given in terms of Gram matrices of periodicity generators offers a good {\em analogy with causal trajectories in  special relativity and Minkowski space-times of arbitrary dimension} \cite{Borel,LEMW}. The analogy works as follows: in a Minkowski space-time, the most conspicuous structural element is the light cone, while in a space of symmetric matrices, the most conspicuous structural element is the positive semidefinite cone; the light cone allows the distinction between causal trajectories and other trajectories in one setting and the positive semidefinite cone allows the distinction of auxetic trajectories from other trajectories in the other setting. In fact, as elaborated below in Section~\ref{sec:Minkowski}, the case $d=2$ is a literal matching.

\begin{corollary}
\label{cor:contractionImpliesDecreasedVolume}
An auxetic path is volume-increasing,  that is, along a (non-trivial) auxetic path, the volume of a fundamental parallelotope (`unit cell') increases. The converse is obviously not true.
\end{corollary}
\noindent
{\em Proof:} \  For a linear contraction operator $T$, we have
 $|det(T)|\leq ||T||^d \leq 1$. An alternative argument uses eigenvalues:  $|det(T)|=\prod_{k=1}^d |\alpha_k|$ is the product of
all eigenvalues (multiplicities included) in absolute value. The volume decrease results from the fact that all eigenvalues satisfy $|\alpha_k| \leq ||T|| \leq 1$ (the spectral radius is bounded by the norm). Thus, for increasing values of the parameter along the auxetic path the volume
goes up. \qed
  
\medskip \noindent
{\bf Remarks.} Structurally, there are additional features which validate our concept of auxetic deformation: (a) the operators  $T_{\tau_2\tau_1}$ are intrinsic, hence the equivalent characterization via Gram matrices does not depend on the choice of periodicity generators,  (b)  for the same reason, upon relaxation of periodicity to a sublattice of finite index $\tilde{\Gamma}\subset \Gamma$, auxetic deformations of the initial framework remain auxetic deformations of the new framework.

\section{Case studies: quartz and cristobalite}
\label{sec:silica}

In this section, we illustrate our mathematical criterion for auxeticity by addressing two fairly classical deformation scenarios: the $\alpha$--$\beta$ phase transition for the silica polymorphs quartz and cristobalite. Our selected references, from an otherwise vast and complex materials science literature \cite{Dol,D,PG}, emphasize the geometric modeling. The simplicity of our method may be compared and contrasted with various other approaches based on wider assumptions and usually dependent on experiment, simulations and multiple numerical approximation algorithms \cite{AE,GW,KC,YHWP}.

\medskip
\noindent
{\bf The $\alpha$--$\beta$ transition of quartz.}
The crystal structure of quartz was determined in the early years of $X$-ray crystallography.
The high temperature phase, called $\beta$--quartz, has higher symmetry and was identified
first. A geometric deformation heuristic was then implicated in the determination of the lower temperature phase, called $\alpha$--quartz \cite{BG,G1}. Later, this deformation, seen as a `coordinated tilting' of the oxygen tetrahedra was explicitly proposed as the underlying geometric mechanism of the $\alpha$--$\beta$ phase transition of quartz \cite{M,GD}.

\medskip
For a direct verification that this particular deformation trajectory is auxetic in the geometric sense described above, it is convenient to start with the $\beta$-configuration. The essentials are illustrated in Fig.~\ref{FigBquartz},  with $\{e_1,e_2, e_3\}$ being the standard orthogonal frame of reference in $R^3$. Our idealized model is made of congruent regular tetrahedra. By periodicity, the three depicted tetrahedra $A, B$ and $C$ are enough for completing the whole (infinite) framework structure of the crystal. The zero sum of the four marked period vectors remains zero in the deformation path to $\alpha$-quartz. This `tilt' scenario rotates the $A_1A_2A_3A_4$ tetrahedron with an angle $\theta$ around the indicated $e_2$ axis and replicates this tilt for the $B$ and $C$ tetrahedra by rotating around the vertical axis $e_3$ with $2\pi/3$ and $4\pi/3$ respectively. Thus, a helical symmetry around vertical directions is maintained. 

\begin{figure}[h]
\centering
 \vspace{-16pt}
 {\includegraphics[width=0.8\textwidth]{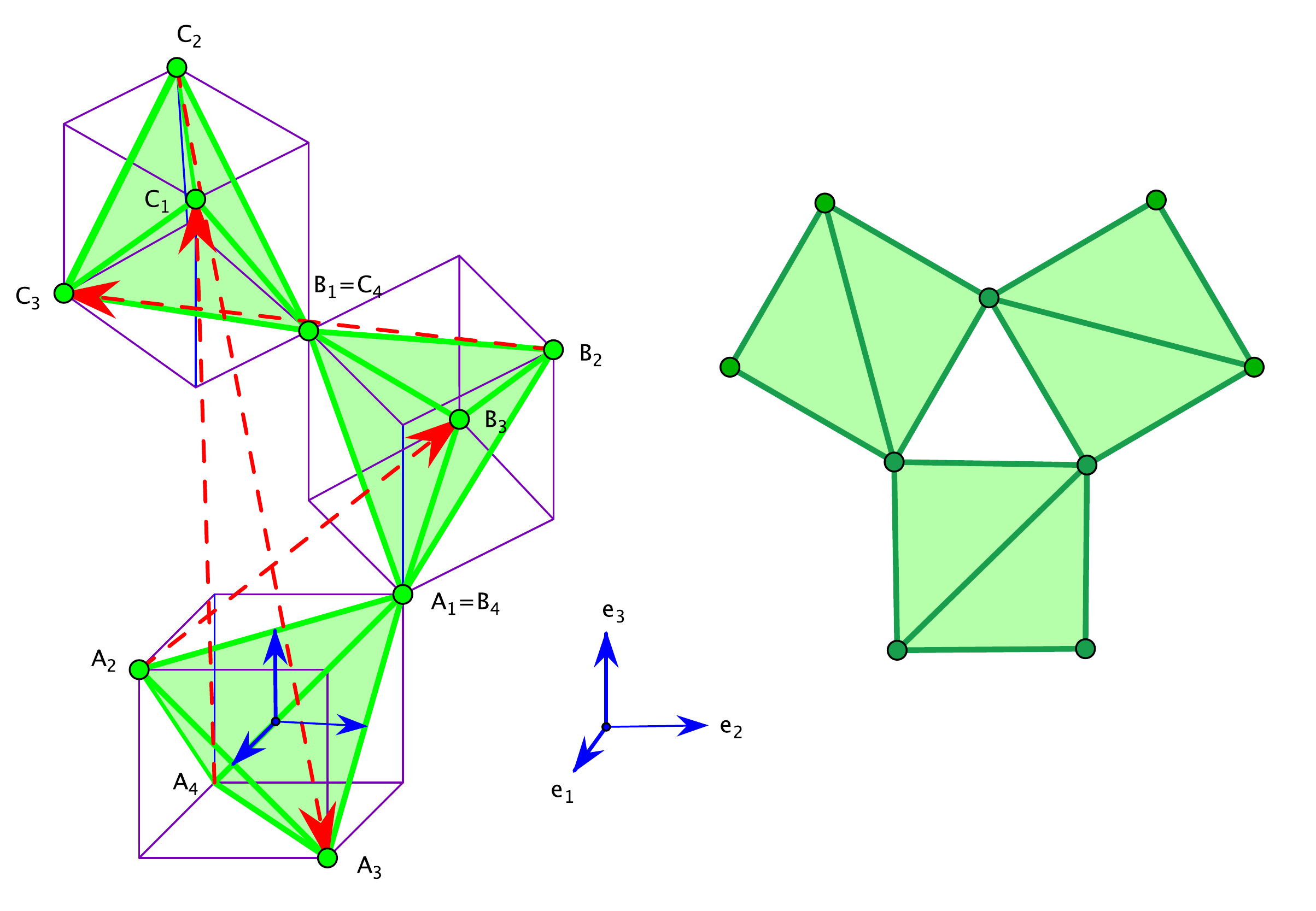}}
 \caption{ A representation of $\beta$--quartz, with highlighted oxygen tetrahedra surrounded by cubes for suggestive purposes. On the right: a view from above.}
 \label{FigBquartz}
 \vspace{-16pt}
\end{figure}

\noindent
We choose $\lambda_1=\overrightarrow{A_2B_3}$, $\lambda_2=\overrightarrow{B_2C_3}$ and $\lambda_3=\overrightarrow{A_4C_1}$ as generators of the periodicity lattice. In order to compute the Gram matrix of this basis, we need the following elements. The `tilting' of $A$ is expressed by a rotation with angle $\theta$ around the second axis. The corresponding matrix is:

\begin{equation}\label{eq:tilt}
T_{\theta}=\left( \begin{array}{ccc} 
\cos{\theta} & 0 & \sin{\theta} \\
0 & 1 & 0 \\
-\sin{\theta} & 0 & \cos{\theta}
\end{array} \right)
\end{equation}

\noindent
$B$ and $C$  are replicas of the tilted $A$ assembled as depicted after rotation around $e_3$ by $2\pi/3$ and $4\pi/3$ respectively. Rotation with $2\pi/3$ gives the matrix:

\begin{equation}\label{eq:rot}
R_{2\pi/3}=\left( \begin{array}{ccc} 
-1/2 & -\sqrt{3}/2 & 0 \\
\sqrt{3}/2 & -1/2 & 0 \\
0 & 0 & 1
\end{array} \right)
\end{equation}

\noindent
By composition we obtain:

\begin{equation}\label{eq:comp}
R_{2\pi/3}T_{\theta}=\left( \begin{array}{ccc} 
-\cos{\theta} /2 & -\sqrt{3}/2 &-\sin{\theta}/2 \\
\sqrt{3}\cos{\theta}/2 & -1/2 & \sqrt{3}\sin{\theta}/2 \\
-\sin{\theta} & 0 & \cos{\theta}
\end{array} \right)
\end{equation}

\noindent
This leads to:

\begin{equation}\label{eq:one}
\overrightarrow{A_2B_3}=\overrightarrow{A_2A_1}+\overrightarrow{B_4B_3}=(1+\sqrt{3}\cos{\theta}) \left( \begin{array}{c} 
-\sqrt{3} \\
1 \\
0
\end{array} \right)
\end{equation}

\noindent
and the resulting Gram matrix has the form:

\begin{equation}\label{eq:Gquartz}
\omega(\theta)=\left( \begin{array}{cc} 
(1+\sqrt{3}\cos{\theta}) ^2 \left( \begin{array}{cc}
4 & -2 \\
-2 & 4
\end{array} \right) & 0 \\
0 & (6\cos{\theta})^2
\end{array} \right)
\end{equation}
 
\noindent
It follows immediately that $\frac{d\omega }{d\theta} (\theta)$ is a negative-definite matrix for $\theta\in (0, \pi/2)$, an interval that includes the $\alpha$ configuration. Thus, the deformation trajectory from $\alpha$ to $\beta$ quartz is auxetic.

\medskip
\noindent
{\bf The $\alpha$--$\beta$ transition of cristobalite.} At the idealized geometric level where we undertake our illustration, the case of cristobalite is quite similar to that of quartz \cite{G2,M,T}. The high temperature phase, also called $\beta$--cristobalite, is presented in Fig.~\ref{FigBcristobalite}. The transition to the low temperature phase, $\alpha$--cristobalite, involves a periodicity lattice which, for the $\beta$--configuration,
is not the maximal lattice of translational symmetries, but an index two sublattice of it. In Fig.~\ref{FigBcristobalite} we mark five period vectors which indicate how to complete the shown fragment to a full crystal structure. These five generators of the periodicity lattice are related by two linear dependence identities respected in the deformation path to $\alpha$--cristobalite. The tilting  of the depicted tetrahedra will maintain a vertical helical symmetry.  This tilt scenario rotates the bottom tetrahedron with an angle $\theta$ around the indicated $e_2$ axis and replicates this tilt for the successive tetrahedra by rotating around the vertical axis $e_3$ with $\pi/2$, $\pi$ and $3\pi/2$ respectively.  

\begin{figure}[h]
\centering
 \vspace{-14pt}
 {\includegraphics[width=0.8\textwidth]{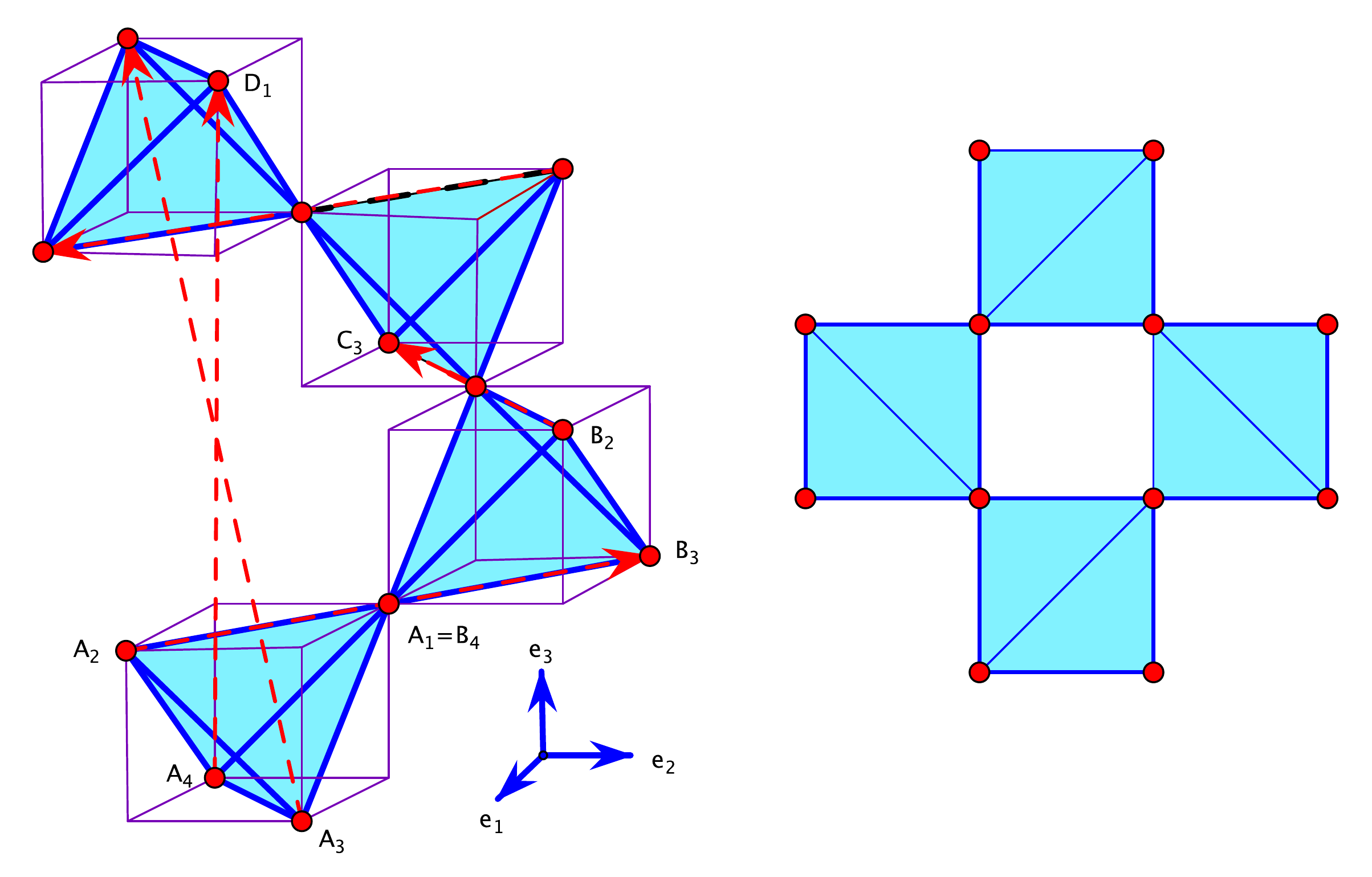}}
 \caption{ A representation of $\beta$--cristobalite, with highlighted oxygen tetrahedra and cubes traced around them for suggestive purposes. On the right: a view from above.}
 \label{FigBcristobalite}
 \vspace{-10pt}
\end{figure} 

\noindent
We choose $\lambda_1=\overrightarrow{A_2B_3}$, $\lambda_2=\overrightarrow{B_2C_3}$ and $\lambda_3=\overrightarrow{A_4D_1}$ as generators of the periodicity lattice. These vectors remain mutually orthogonal as the periodicity lattice varies with the deformation. A computation  similar to the one conducted above for quartz yields the following Gram matrix.

\begin{equation}\label{eq:Gcristobalite}
\omega(\theta)=\left( \begin{array}{ccc} 
8(1+\cos{\theta})^2 & 0 & 0 \\
0 & 8(1+\cos{\theta})^2 & 0 \\
0 & 0 & (8\cos{\theta})^2
\end{array} \right)
\end{equation}

\noindent
The derivative $\frac{d\omega }{d\theta} (\theta)$ is obviously negative-definite for $\theta \in (0,\pi/2$, hence the $\alpha$ to $\beta$ deformation trajectory of cristobalite is auxetic.

\section{Planar frameworks: expansive and auxetic}
\label{sec:planar}

This section explores the rapport between expansive and auxetic deformations in dimension two, where it has far-reaching consequences for auxetic design. Expansiveness is the stronger property and, in dimension two, it is intimately related to a remarkable class of periodic frameworks called pointed pseudo-triangulations. We show that an elementary procedure for generating periodic pointed pseudo-triangulations leads to an endless series of expansive and thereby auxetic designs. This alters utterly the catalog of planar auxetic structures which, up to this point, was sparsely populated with a few dozen examples\cite{EL}. In fact, several well-known examples in the older catalog are easily recognized as frameworks with auxetic capabilities via refinements to pseudo-triangulations.
 
\medskip 
\noindent
{\bf Expansive one-parameter deformations.}\ A one-parameter deformation of a periodic framework is called {\em expansive} when all the distances between pairs of vertices increase or stay the same (when the parameter increases). Thus, it is not possible for one pair of vertices to get closer together while another pair moves further apart.

\begin{figure}[h]
\centering
\subfloat[Pointed]
 {\includegraphics[width=0.23\textwidth]{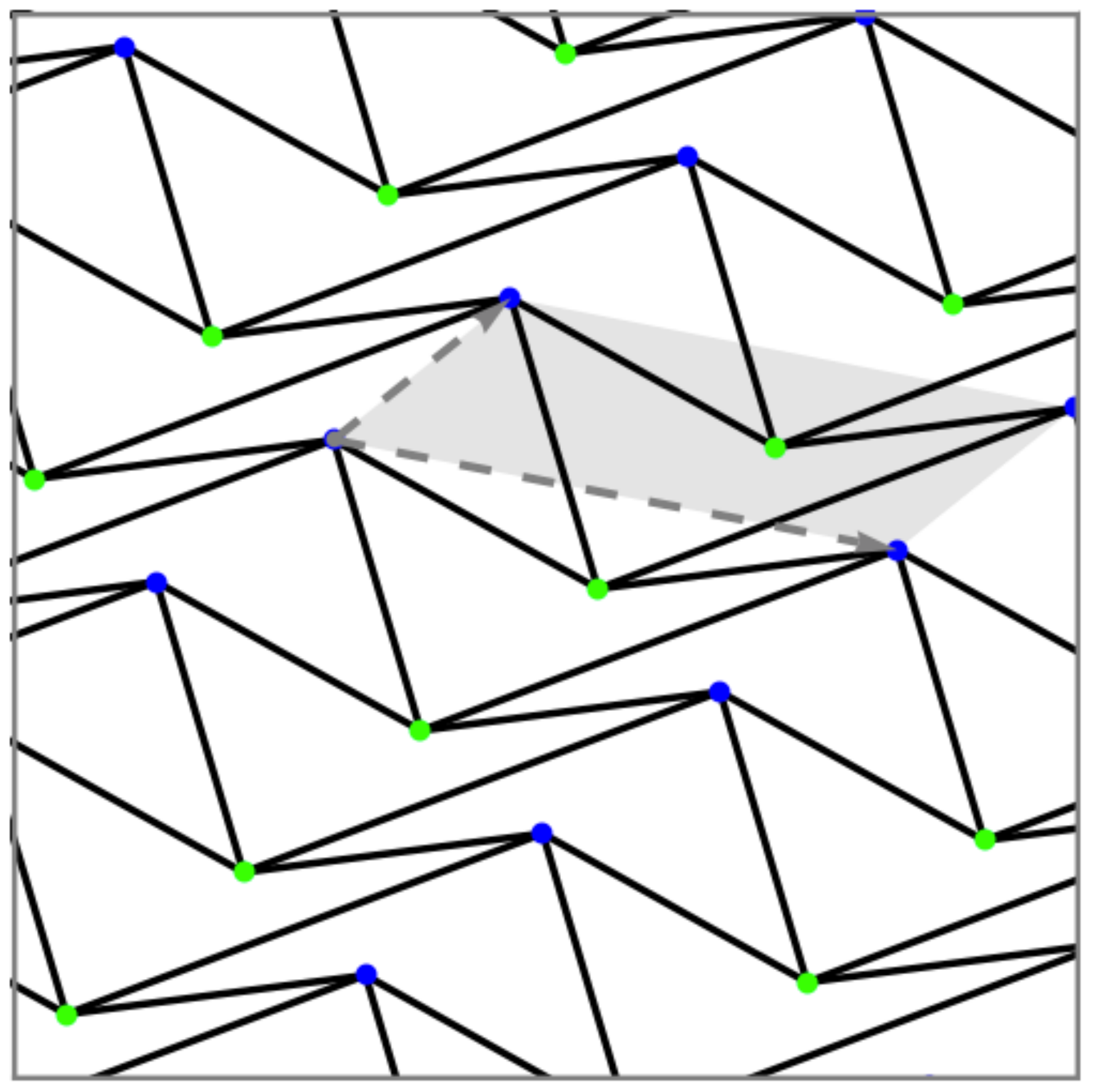}}
 \hspace{3pt}
 \subfloat[Pointed] {\includegraphics[width=0.23\textwidth]{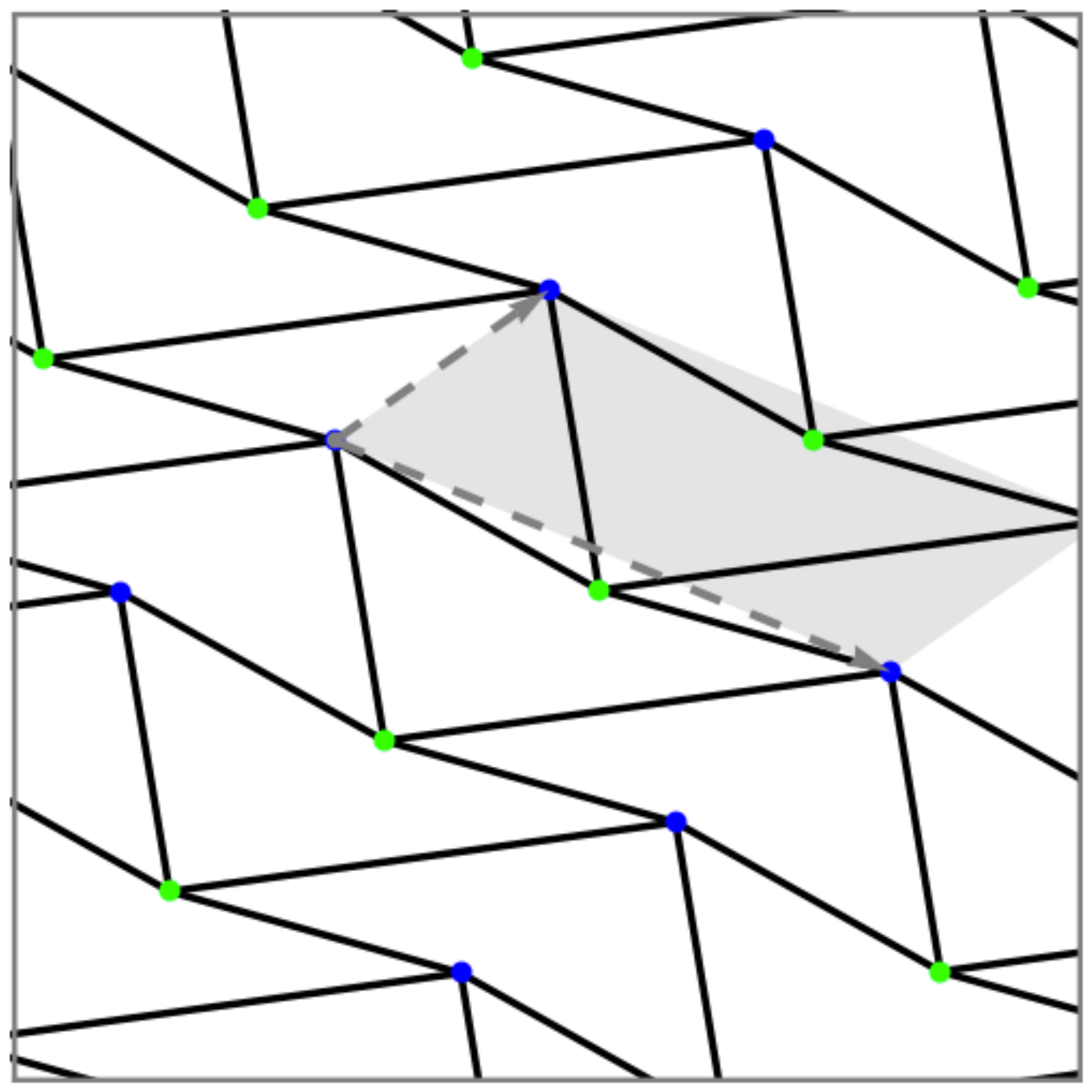}}
 \hspace{3pt}
 \subfloat[Alignment] {\includegraphics[width=0.23\textwidth]{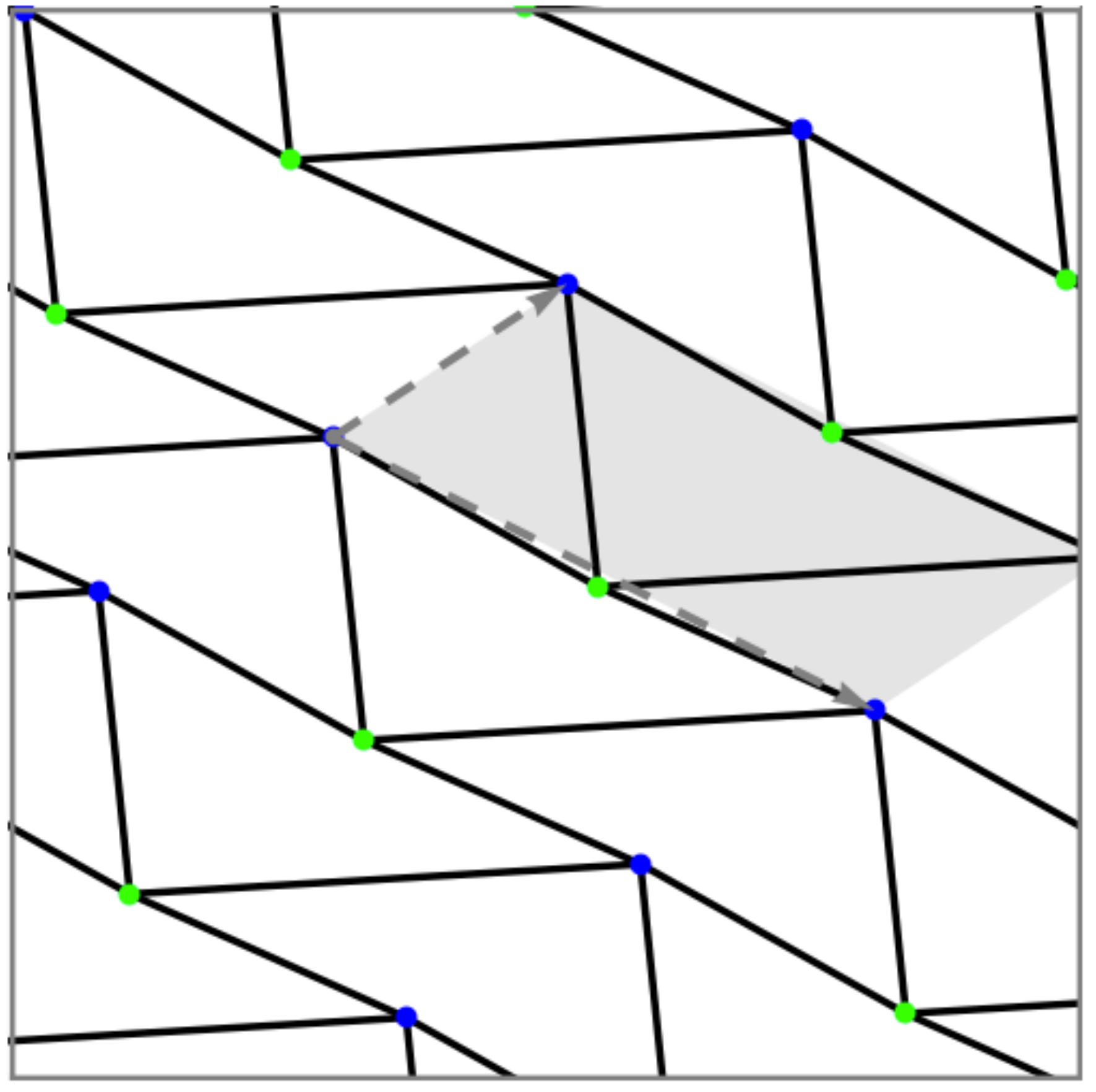}}
 \hspace{3pt}
 \subfloat[Not pointed] {\includegraphics[width=0.23\textwidth]{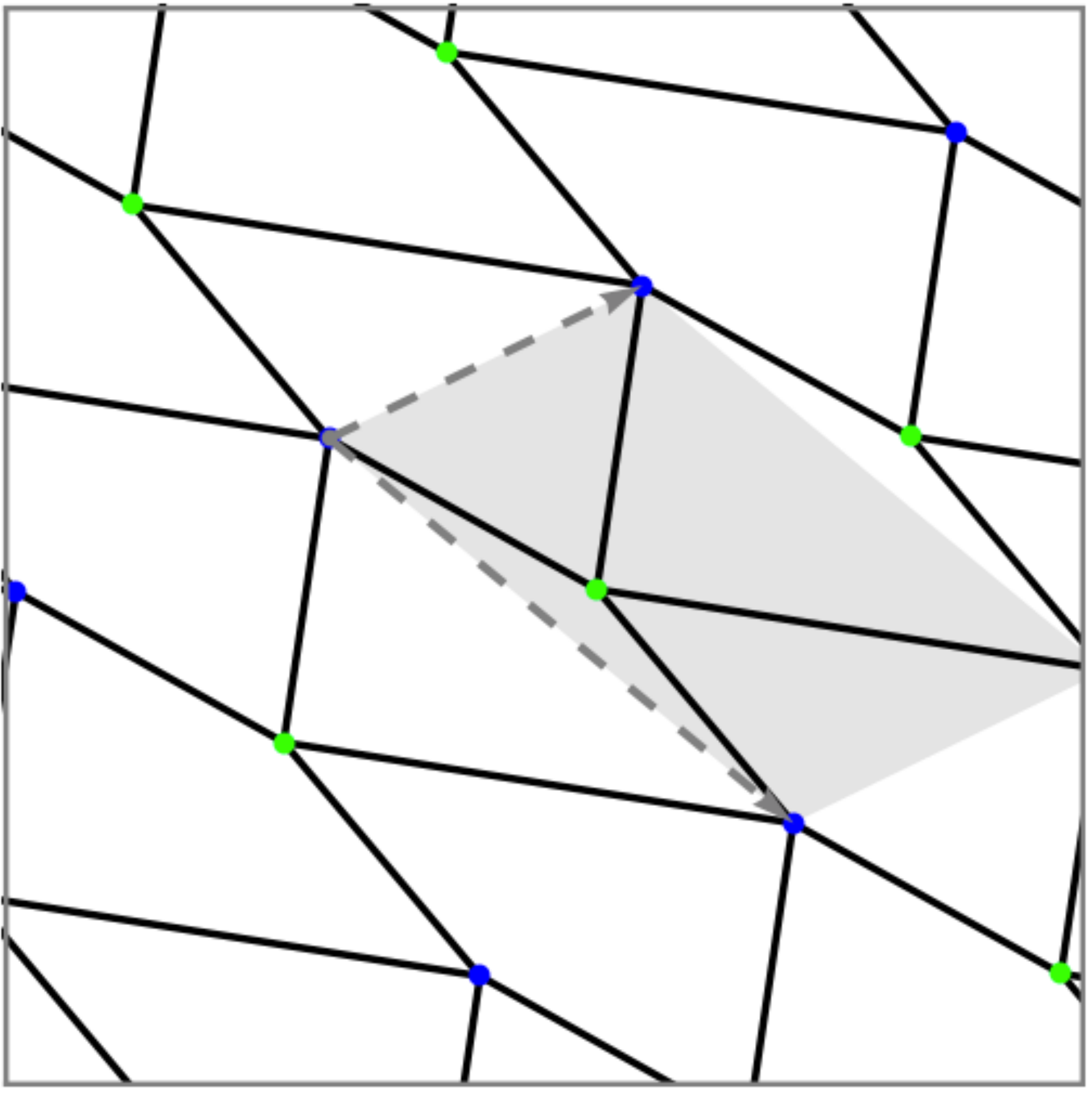}}
 \caption{(a,b) Two snapshots of the unique, expansive trajectory of a periodic pseudo-triangulation. When the mechanism ceases to be pointed through the alignment of edges (c), both the expansive and the auxetic behavior disappear. (d) In this configuration, the lattice generators have opposite growth behavior: an increase along one axis leads to a decrease along the other axis.} 
 \label{FigAnimPPT}
\end{figure}

\medskip \noindent
{\bf Periodic pseudo-triangulations.}\ A {\em pseudo-triangle} is a simple closed planar polygon with exactly three internal angles smaller than $\pi$. A set of vectors (without aligned vectors in opposite directions)  is {\em pointed} if there is no linear combination with strictly positive coefficients that sums them to $0$.  Equivalently, a set of vectors with the same origin is pointed when contained in some open half-plane determined by a line through their common origin. Thus, for a pointed set of vectors, some consecutive pair (in the circular rotational order around the common origin) has an angle larger than $\pi$.  Examples of pointed and non-pointed periodic frameworks appear in Figs.~\ref{FigAnimPPT} and \ref{FigReentrantH}. A planar non-crossing periodic framework is a {\em periodic pointed pseudo-triangulation} (shortly, a {\em periodic pseudo-triangulation}) when all faces are pseudo-triangles and the framework is pointed at every vertex.  Such periodic frameworks are {\em maximal} with the property of being non-crossing, pointed, flexible and non-redundant (independent) (in the sense of rigidity theory, see \cite{S2, BS4}), and any added edge leads to a violation of one or more of these properties. See Figs.~\ref{FigAnimPPT}, \ref{FigReentrantH} and \ref{FigPPT}.

\medskip
\noindent
{\bf Periodic pseudo-triangulations are expansive mechanisms.}  The most remarkable property  of a periodic pointed pseudo-triangulation (proven in \cite{BS4}) is that it is flexible with exactly one-degree-of-freedom and has an expansive one-parameter deformation for as long as the deformed framework remains a pseudo-triangulation. An illustration of this phenomenon is shown in Fig.~\ref{FigAnimPPT}. A periodic pseudo-triangulation also has the remarkable property of remaining a one-degree-of-freedom expansive mechanism for {\em any relaxation of periodicity.} Furthermore, any infinitesimal expansive deformation of a pointed non-crossing framework is obtained as a convex combination of its refinements to pseudo-triangulations. The proof of these facts relies on several advanced techniques, including a generalization (from finite to periodic) of the classical theorem of James Clerk Maxwell \cite{Max} concerning stresses and liftings of planar frameworks.

\medskip
\noindent
{\bf Expansive implies auxetic.}\ The fact that {\em expansive implies auxetic} in our geometric sense is proven with a short argument.

\begin{theorem}\label{expaux}
Let $(G,\Gamma, p_{\tau},\pi_{\tau}),\tau\in (-\epsilon,\epsilon)$ be a one-parameter
deformation of a periodic framework in $R^d$. If the path is expansive, that is, if the
 distance between any pair of vertices increases or stays the same for increasing $\tau$, then
the path is also auxetic. However, auxetic paths need not be expansive.
\end{theorem}

\noindent
{\em Proof:}\ The auxetic property depends only on the curve $\omega(\tau)$ and it will be enough to use the expansive property on one orbit of vertices. We have to verify that the operator $T_{\tau_2\tau_1}$ which takes the period lattice basis $\Lambda_{\tau_2}$ to the period lattice basis $\Lambda_{\tau_1}$ is a contraction for $\tau_2 > \tau_1$.

\medskip 
\noindent
In the unit ball of $R^d$, the vectors with rational coordinates relative to the basis $\Lambda_{\tau_2}$ give a dense subset. Since some integer multiple of such a point is a period at moment $\tau_2$, and this period, as a distance between two vertices in a vertex orbit, can only decrease or preserve its norm when mapped by $T_{\tau_2\tau_1}$ to the corresponding period at moment $\tau_1$, we see that a dense subset of points in the unit ball must be mapped to the unit ball. This is enough to conclude that $||T||\leq 1$. \qed
 
\begin{figure}[h]
\centering
 \subfloat[]{\includegraphics[width=0.31\textwidth]{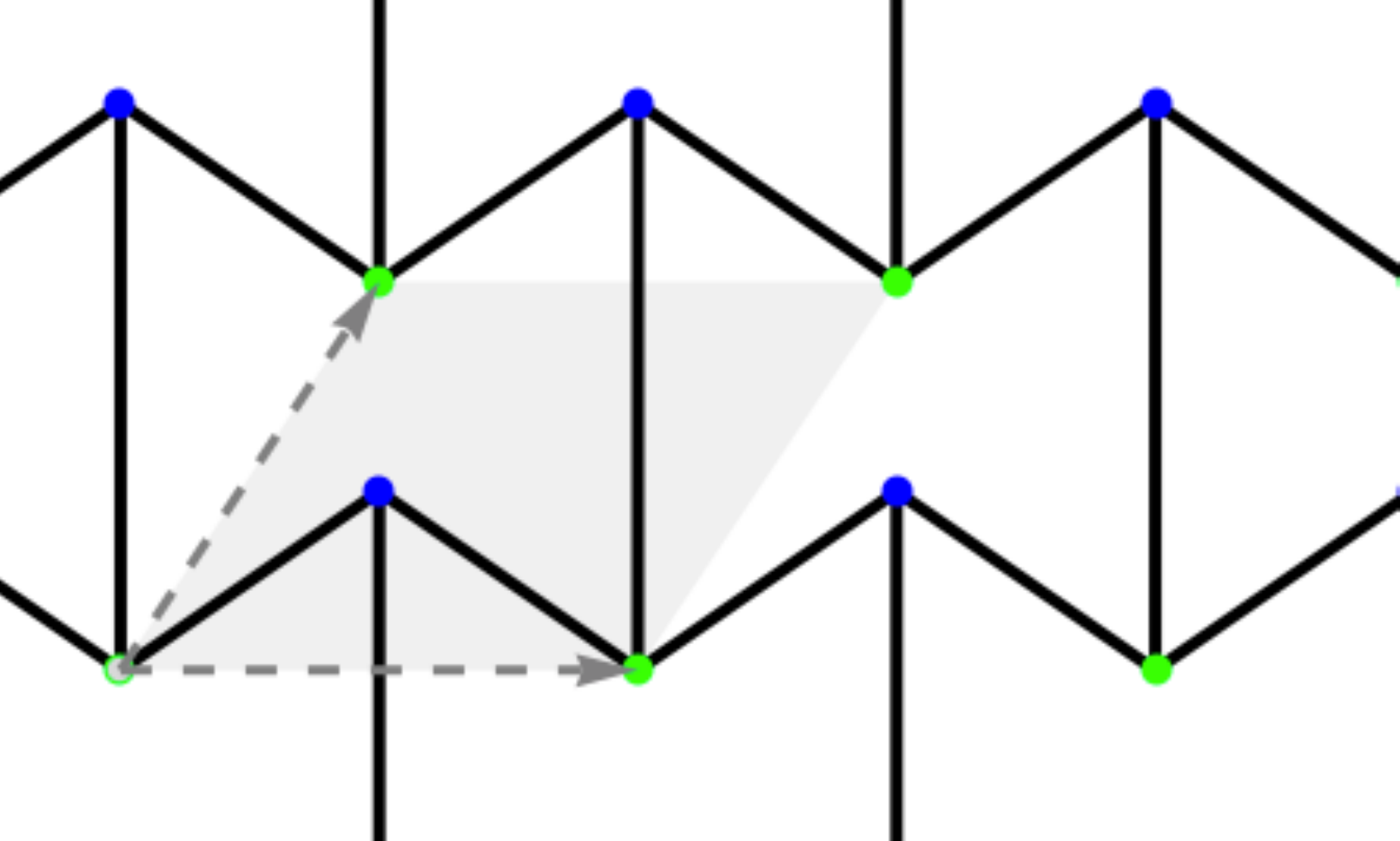}}\hspace{5pt}
 \subfloat[]{\includegraphics[width=0.31\textwidth]{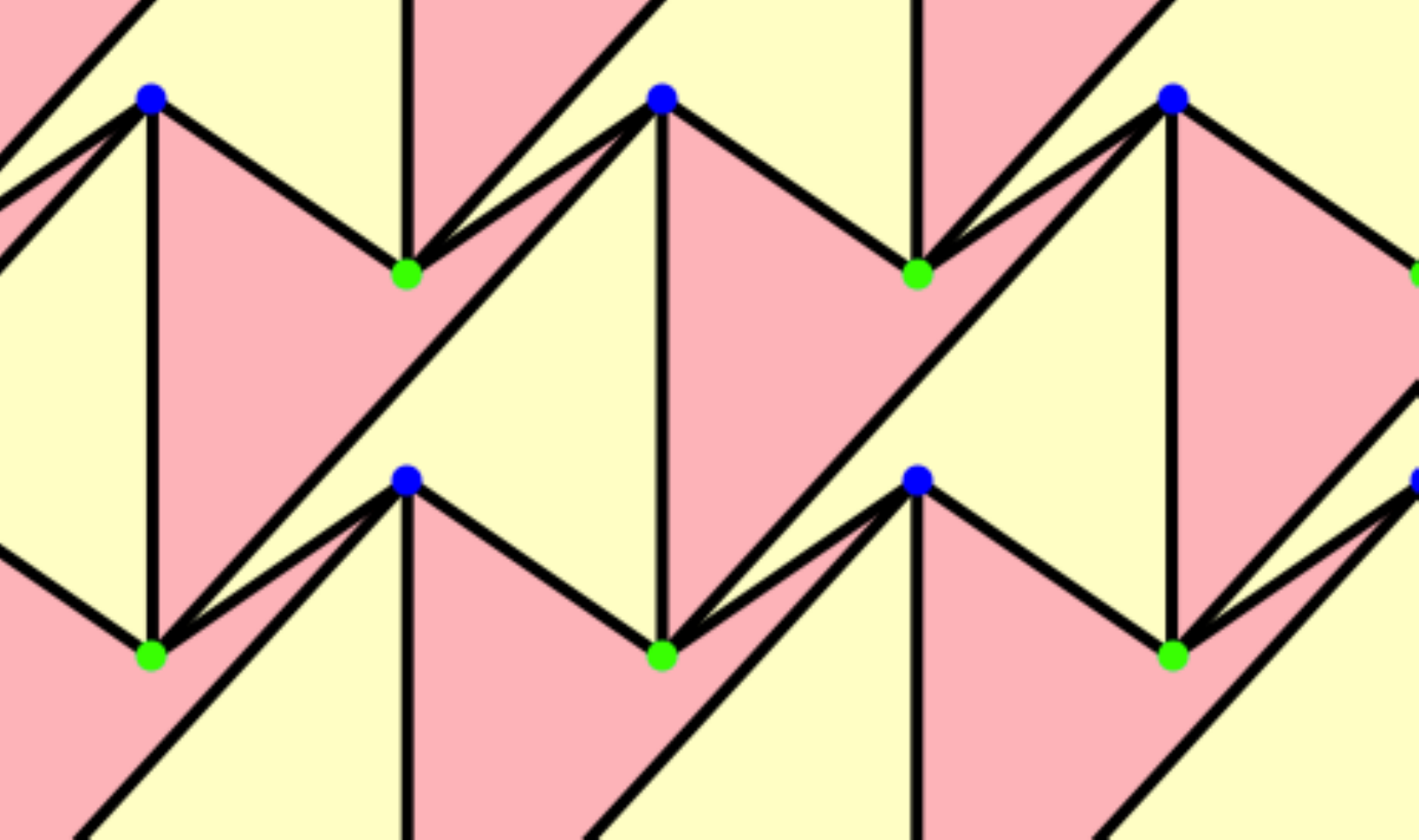}}\hspace{5pt}
 \subfloat[]{\includegraphics[width=0.31\textwidth]{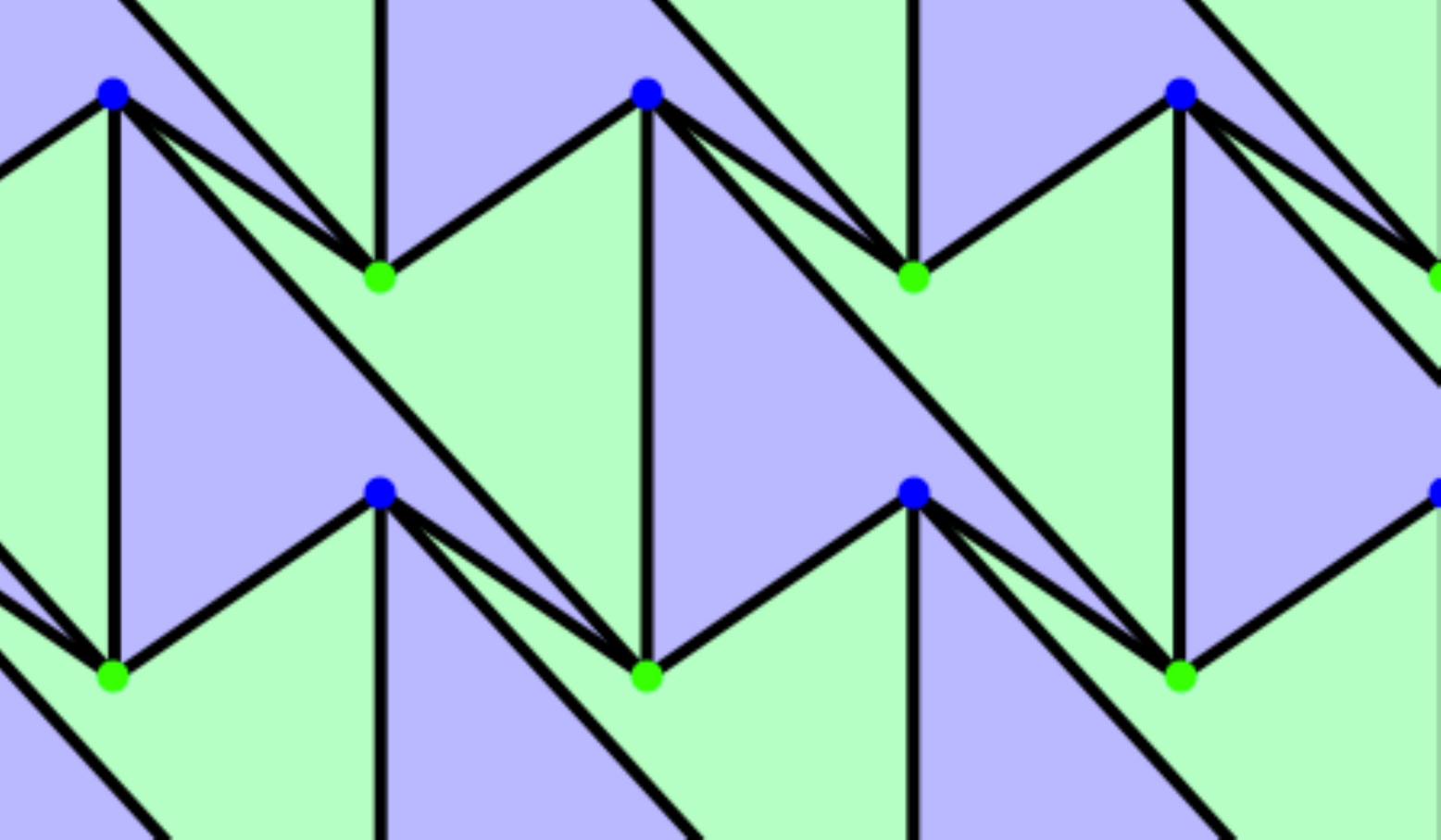}}\\
 \vspace{-4pt}
 \subfloat[]{\includegraphics[width=0.31\textwidth]{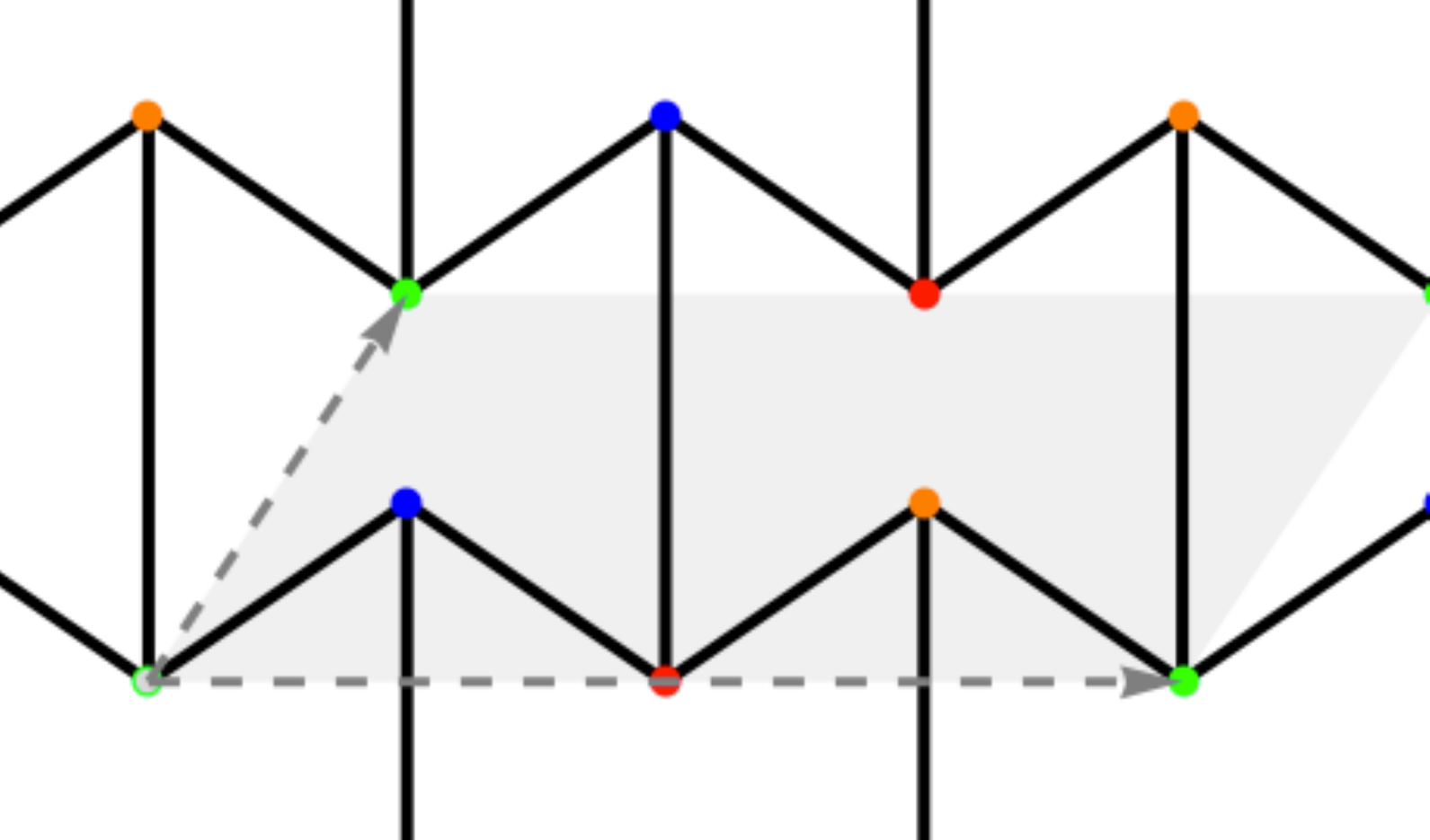}}\hspace{5pt}
 \subfloat[]{\includegraphics[width=0.31\textwidth]{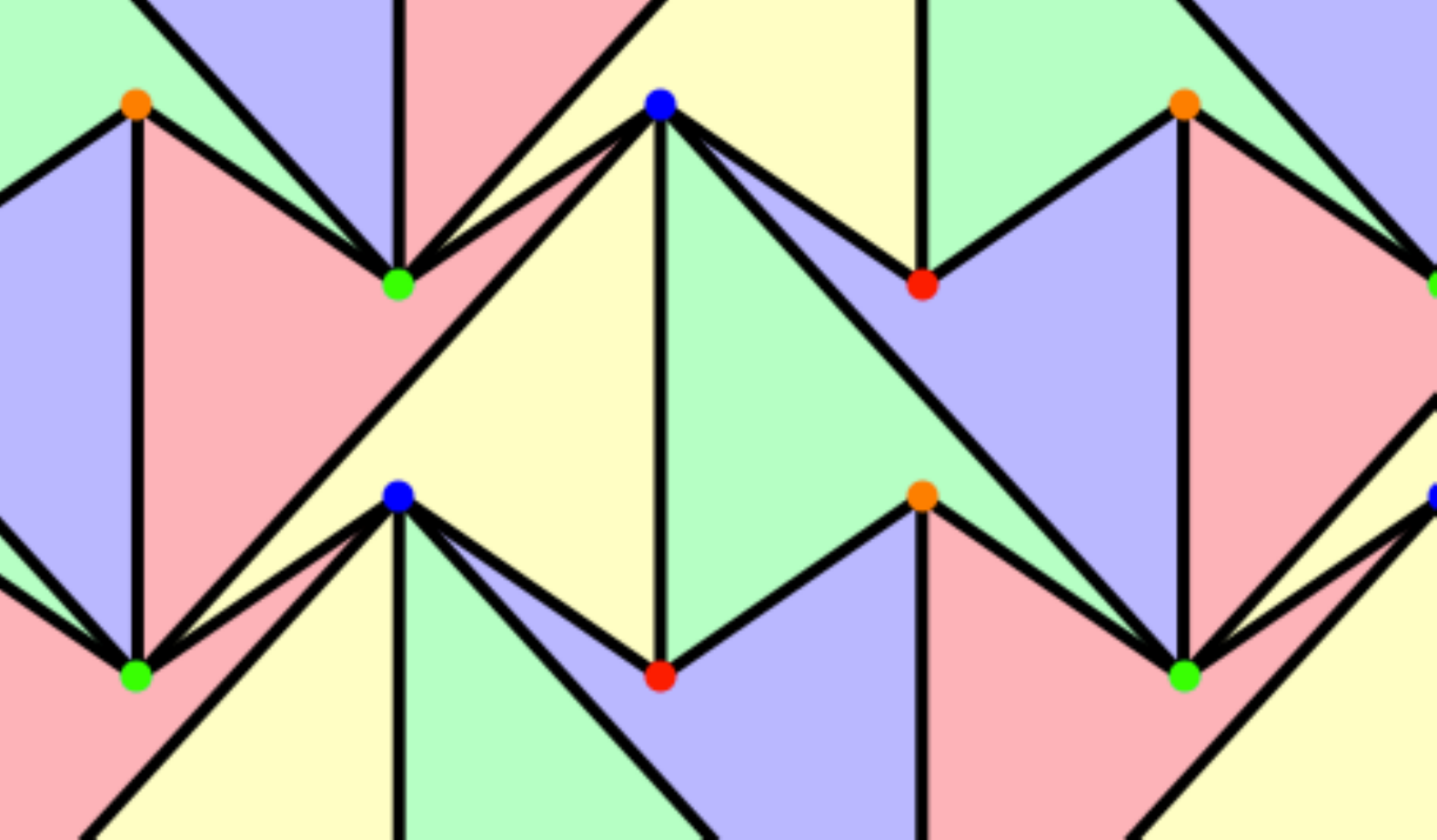}}\hspace{5pt}
 \subfloat[]{\includegraphics[width=0.31\textwidth]{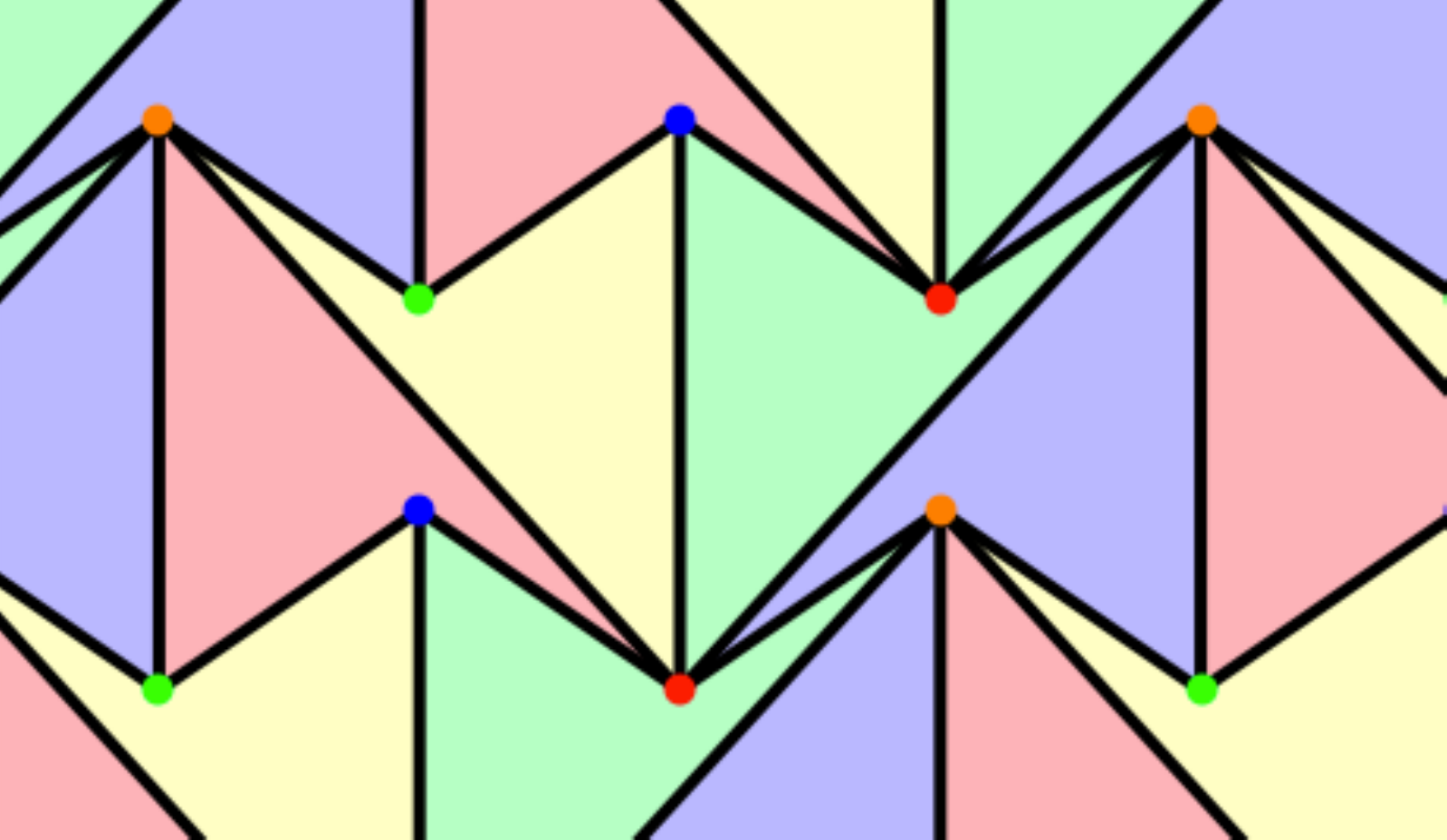}}\\
 \caption{ (a) A hexagonal pointed honeycomb, taken with maximal periodicity and (b, c) its two possible refinements to periodic pseudo-triangulations. Its auxetic capabilities follow from the expansive properties inherent in pseudo-triangulations. (d) A relaxation of the lattice on the same infinite graph increases the degrees-of-freedom from two to three and the new periodic framework has two additional pseudo-triangular refinements (e,f).}
 \label{FigReentrantH}
\end{figure}

\medskip  The Kagome type of framework offers simple examples of configurations which remain auxetic  after ceasing to be expansive \cite{BS4} (pg. 777). Proposition~\ref{prop:hexagonal} in Section \ref{sec:Minkowski} below will show the distinction of these two notions for hexagonal honeycomb frameworks.

\medskip
\noindent
{\bf Hexagonal honeycombs.} In Fig.~\ref{FigReentrantH} we illustrate one of the iconic examples of planar periodic structures with auxetic behavior, the so-called ``reentrant honeycomb'' \cite{Alm,ENHR,Kol}. Our geometric theory of periodic pseudo-triangulations and auxetic deformations shows that the structural source of its expansive and auxetic capabilities resides in the possibility of refining this pointed structure to periodic pseudo-triangulations, by inserting new orbits of edges. For maximal periodicity (Fig.~\ref{FigReentrantH}(a,b,c)), there are two distinct ways for completion to a periodic pseudo-triangulation, and each one induces an expansive trajectory.  Since the honeycomb has two degrees-of-freedom, its deformation space is two-dimensional (a surface), and each expansive trajectory induced by a pseudo-triangulation is a curve on this surface. At the infinitesimal level, i.e. in the two-dimensional tangent plane of the deformation surface, the tangents to the two curves are the two extremal rays of a cone, called the ``expansive cone'': any infinitesimal expansive motion lies in this cone. 

\medskip
\noindent
Fig.~\ref{FigReentrantH}(d,e,f) illustrates the role of the periodicity group $\Gamma$ in the analysis of degrees-of-freedom and structure of the cones of expansive and auxetic infinitesimal deformations of a periodic framework. A relaxation of the maximal periodicity group $\Gamma$ in (a) to an index $2$ sublattice of periods in (d) allows new deformations. Initially we have $n=2$, $m=3$ and $f=2n+1-m=2$ degrees of freedom, while after relaxation we have $\tilde{n}=2n=4$, $\tilde{m}=2m=6$ and $\tilde{f}=2\tilde{n}+1-\tilde{m}=3$ degrees-of-freedom. The cone of expansive infinitesimal deformations has now four extremal rays, corresponding to the indicated refinements to pseudo-triangulations.

\medskip
\noindent
{\bf `Missing rib' models.} The example in Fig.~\ref{FigMissingRib}(a) is related to a so-called `missing rib' planar framework akin to the foam structure considered in \cite{GRSGE}, sample 3.  The figure shows the four new edge orbits which have to be inserted for obtaining a framework {\em kinematically equivalent} (in the sense of \cite{BS6}) with the one depicted in Fig.~\ref{FigPseudoTr}(b). This means that substructures that are rigid may fail to be pointed as long as they are convex, without changing the overall expansive properties of a framework. Indeed, we simply replace them by pointed triangulated convex polygons. In going from Fig.~\ref{FigPseudoTr}(a) to (b), we replaced the connecting squares (having an internal vertex) with an equivalent, rigid, triangulated square.

\begin{figure}[h]
\centering
\subfloat[] {\includegraphics[width=0.44\textwidth]{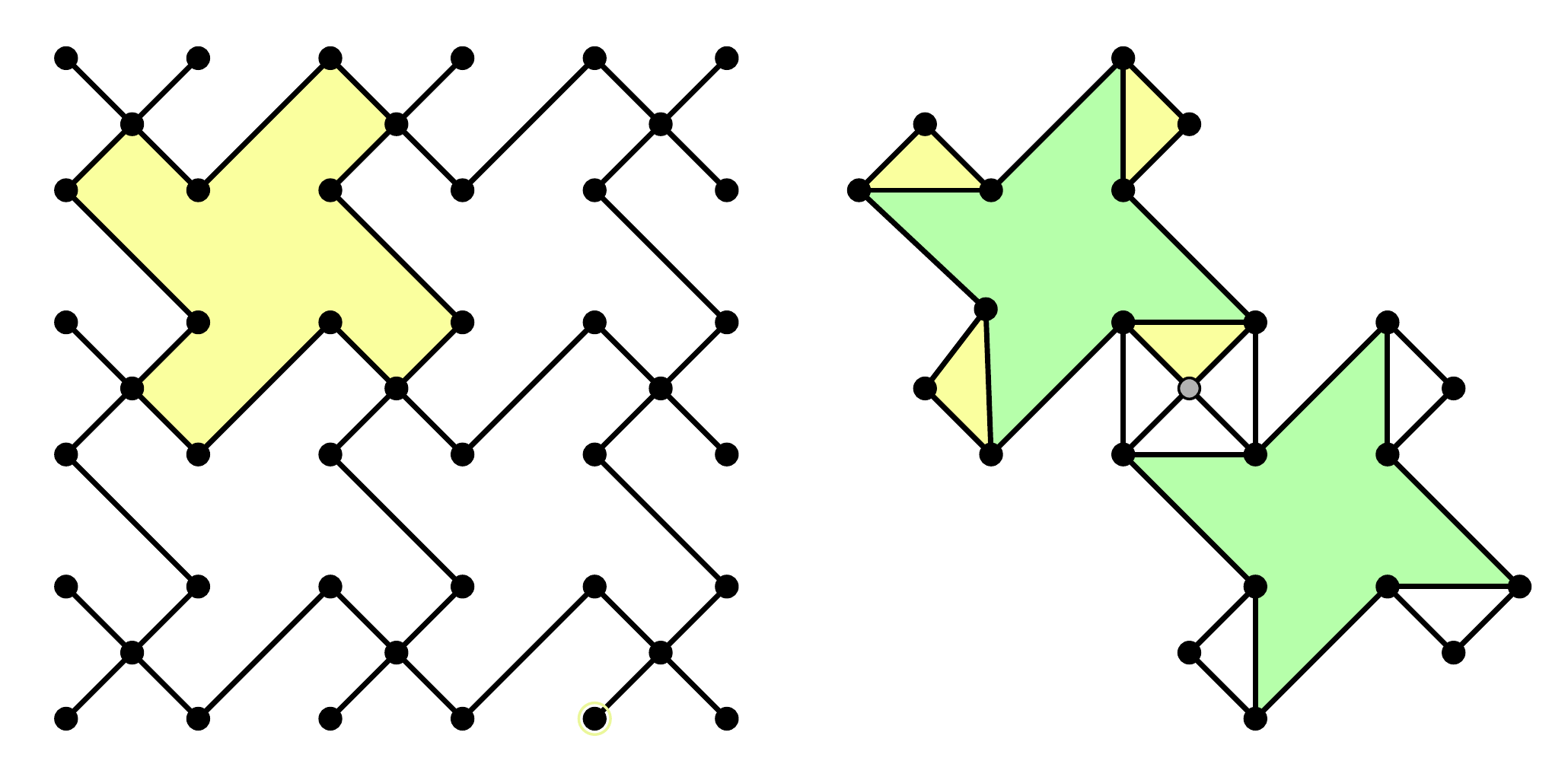}}
 \subfloat[]{\includegraphics[width=0.55\textwidth]{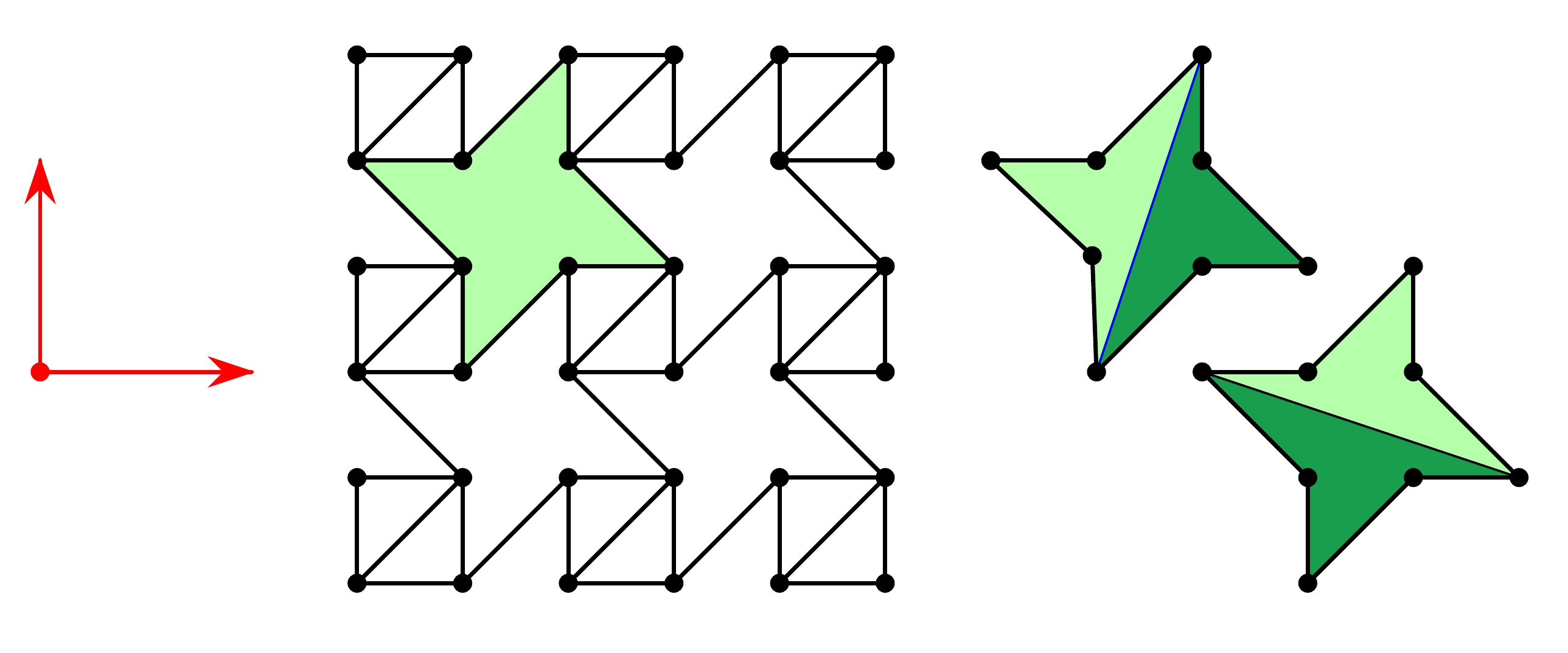}}
 \caption{(a) A `missing rib' model as discussed in \cite{GRSGE} and its conversion to the framework shown nearby.
 (b) A planar periodic framework with two distinct refinements to pseudo-triangulations.
 }
 \label{FigMissingRib}
 \label{FigPseudoTr}
\end{figure}

\medskip
\noindent
With the indicated generators of the periodicity lattice, the framework in Fig.~\ref{FigPseudoTr}(b) has $n=4$ vertex orbits and $m=7$ edge orbits. From the dimension counts proven in \cite{BS1}, it follows that this framework has two degrees of freedom, more precisely a smooth two-dimensional local deformation space. Similar considerations give five degrees of freedom for the `missing rib' framework in Fig.~\ref{FigMissingRib}(a).  Because the second framework is obtained by adding rigid bars (edges) to the first one, it follows that the two-dimensional deformation space of (b) is included in the five-dimensional deformation space of (a).

\medskip
\noindent
As in the reentrant honeycomb case, the framework in Fig.~\ref{FigPseudoTr}(b) has expansive capabilities: indeed, since it is pointed, refinements to pseudo-triangulations are possible, as shown. The `missing rib' structure clearly inherits all these expansive trajectories. {\em This is a fast way to confirm auxetic capabilities in the latter structure, without any explicit calculations.} In \cite{GRSGE}, a hint at possible auxetic behavior is obtained after a considerable more complex exploration, mixing experiment, simulation and (sometimes conflicting) computational results. The role of multiple degrees of freedom is not clarified in the cited paper.

\medskip
\noindent
{\bf Other examples.}  The auxetic property of the ``tetramer system"  considered in \cite{Tr}  can be explained by an underlying framework similar to the framework in Fig.~\ref{FigPseudoTr}(b).

\medskip
\noindent
{\bf Generating expansive trajectories.} Computationally,  a deformation trajectory is generated by numerical integration via calculation of infinitesimally motions at any time step.  The previous examples can be generalized to the following procedure for calculating {\em arbitrary} expansive (and thus, auxetic) trajectories for any  pointed and non-crossing periodic framework, not just pseudo-triangulations. When the framework has more than one degree-of-freedom (as is the case with the ``reentrant honeycomb''), refinements to pseudo-triangulations induce the extremal rays of the cone of infinitesimal expansive motions; an arbitrary infinitesimal expansive motion is then computed as a convex combination of those corresponding to the extremal rays.  

\medskip
The method can sometimes also be applied to situations when the framework has edge crossings or is not pointed (as illustrated in the ``missing rib'' example above). In this case, we seek to turn the substructure violating pointedness or non-crossing into a convex rigid unit (called {\em rigid component} in the rigidity theory literature). We then replace it with a kinematically equivalent pointed pseudo-triangulation, as illustrated in going from Fig.~\ref{FigPseudoTr}(a) to Fig.~\ref{FigPseudoTr}(b). Further details appear in \cite{BS6}.

\begin{figure}[h]
\centering
 {\includegraphics[width=0.23\textwidth]{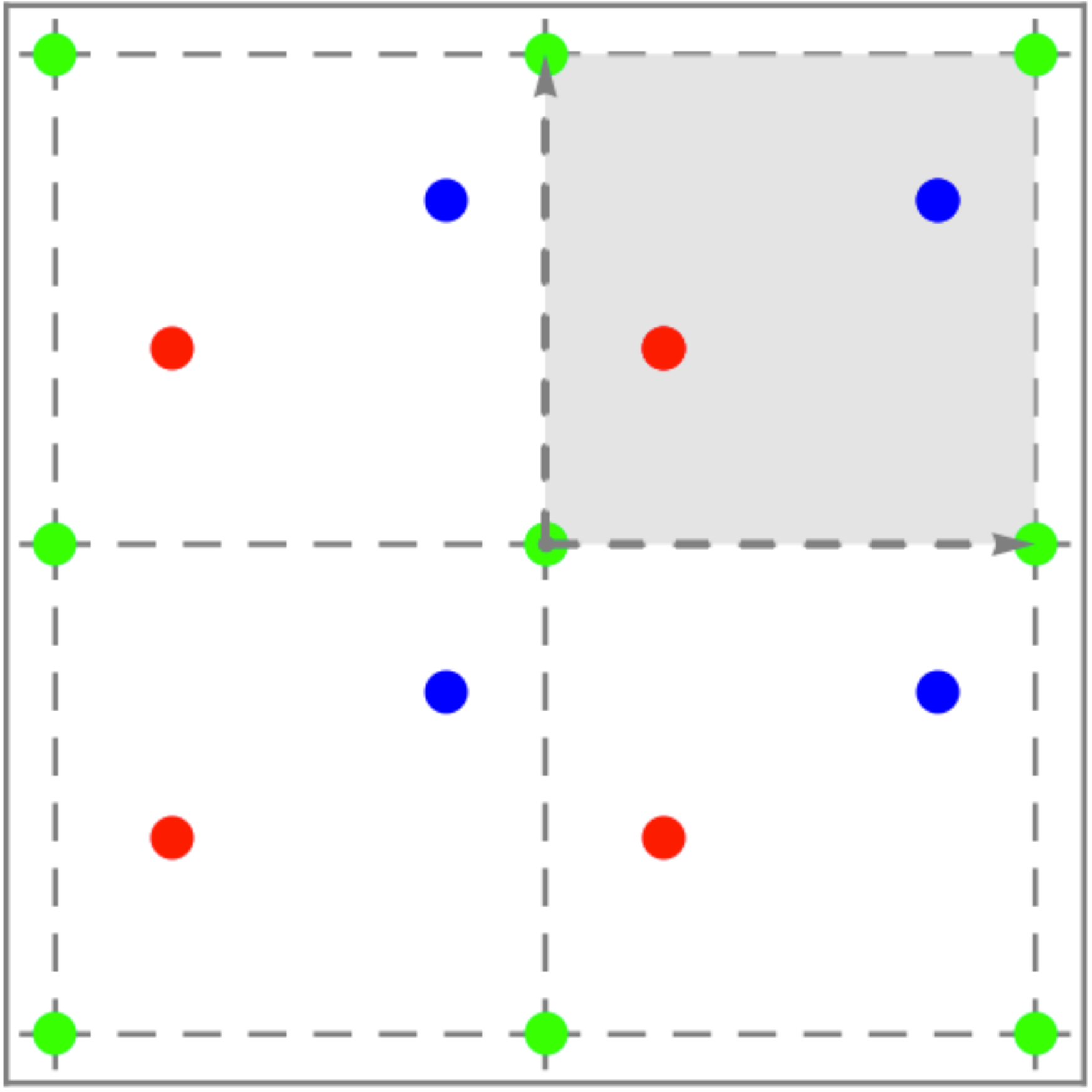}}\hspace{3pt}
 {\includegraphics[width=0.23\textwidth]{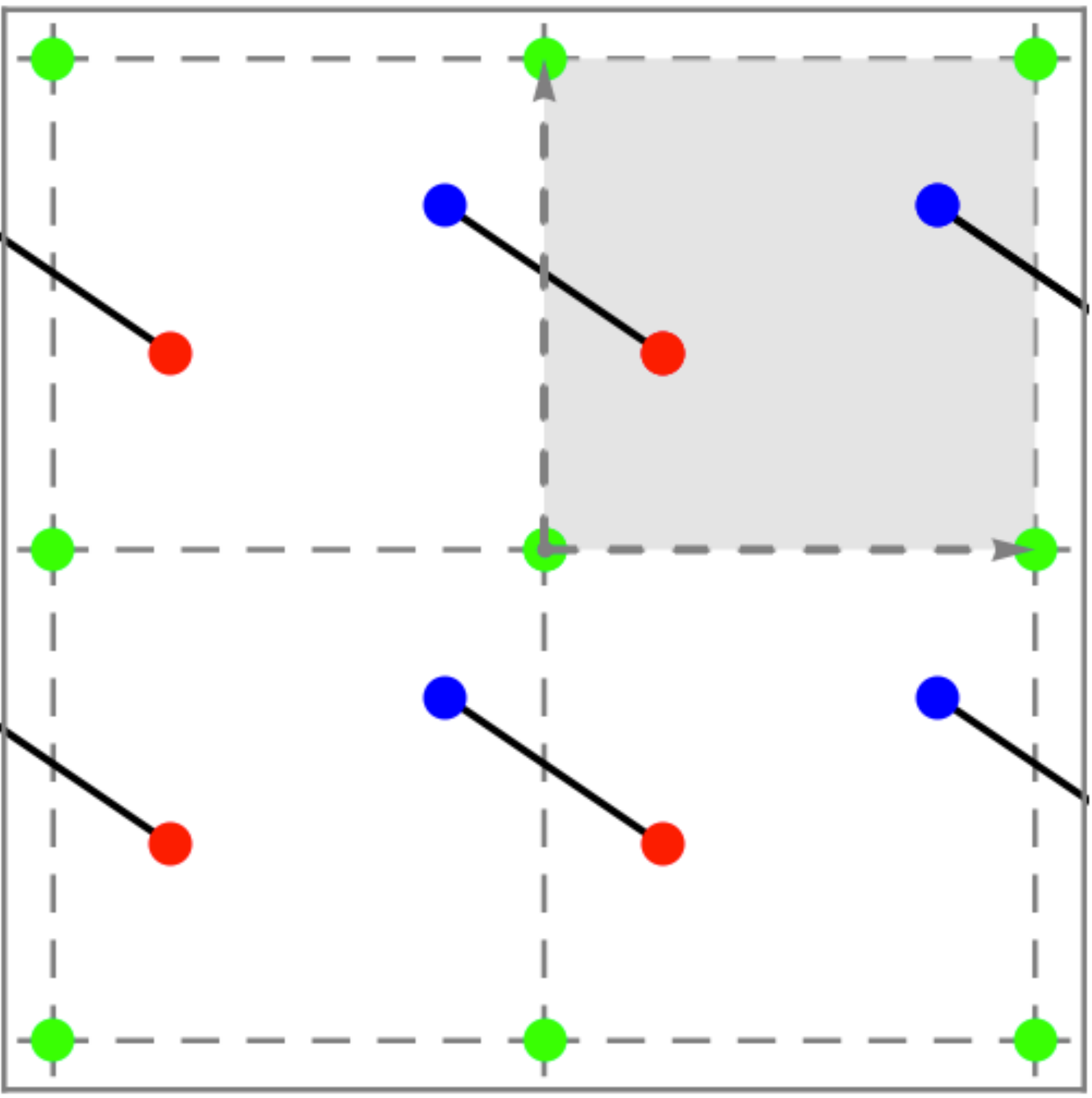}}\hspace{3pt}
 {\includegraphics[width=0.23\textwidth]{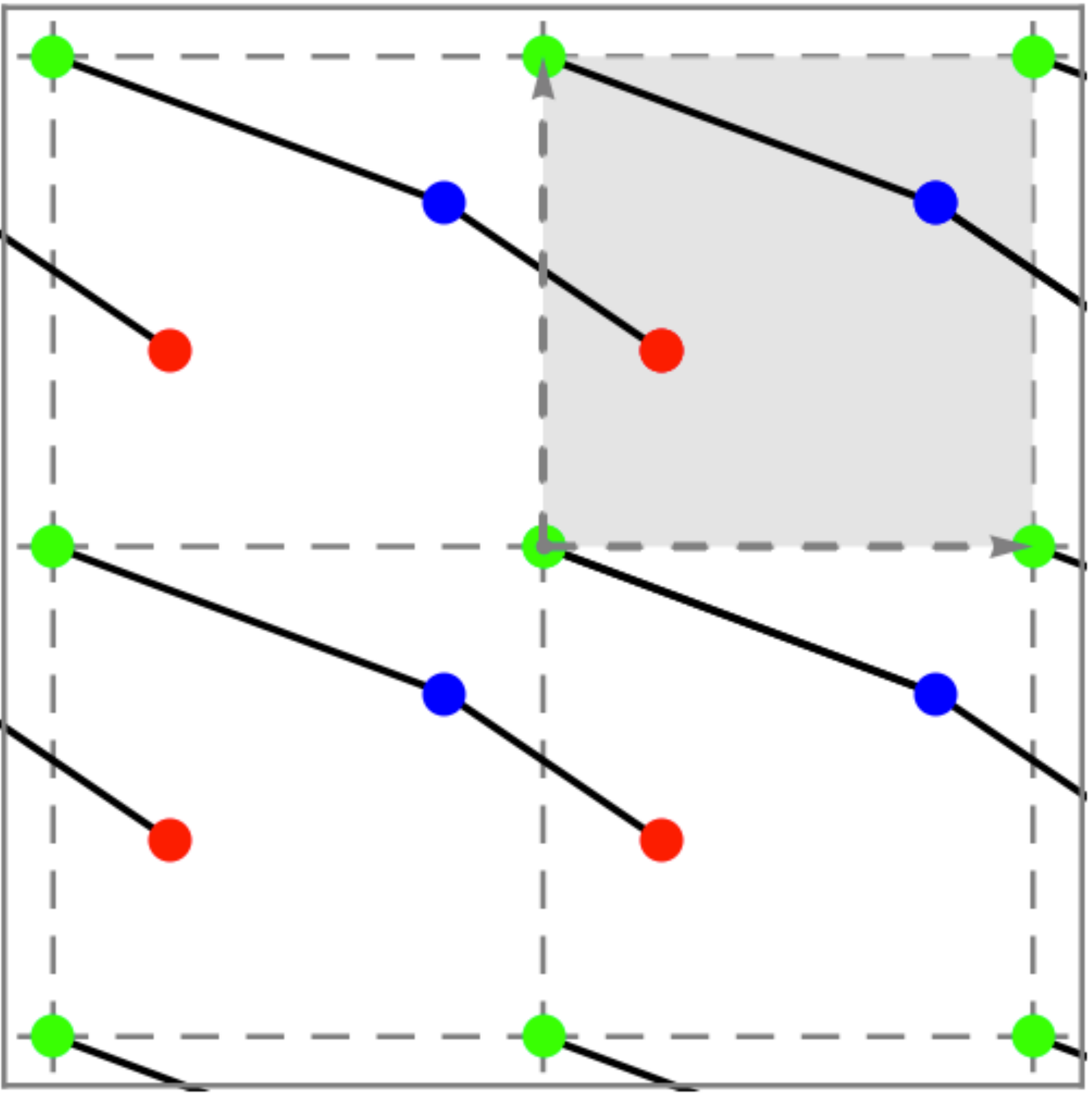}}\hspace{3pt}
 {\includegraphics[width=0.23\textwidth]{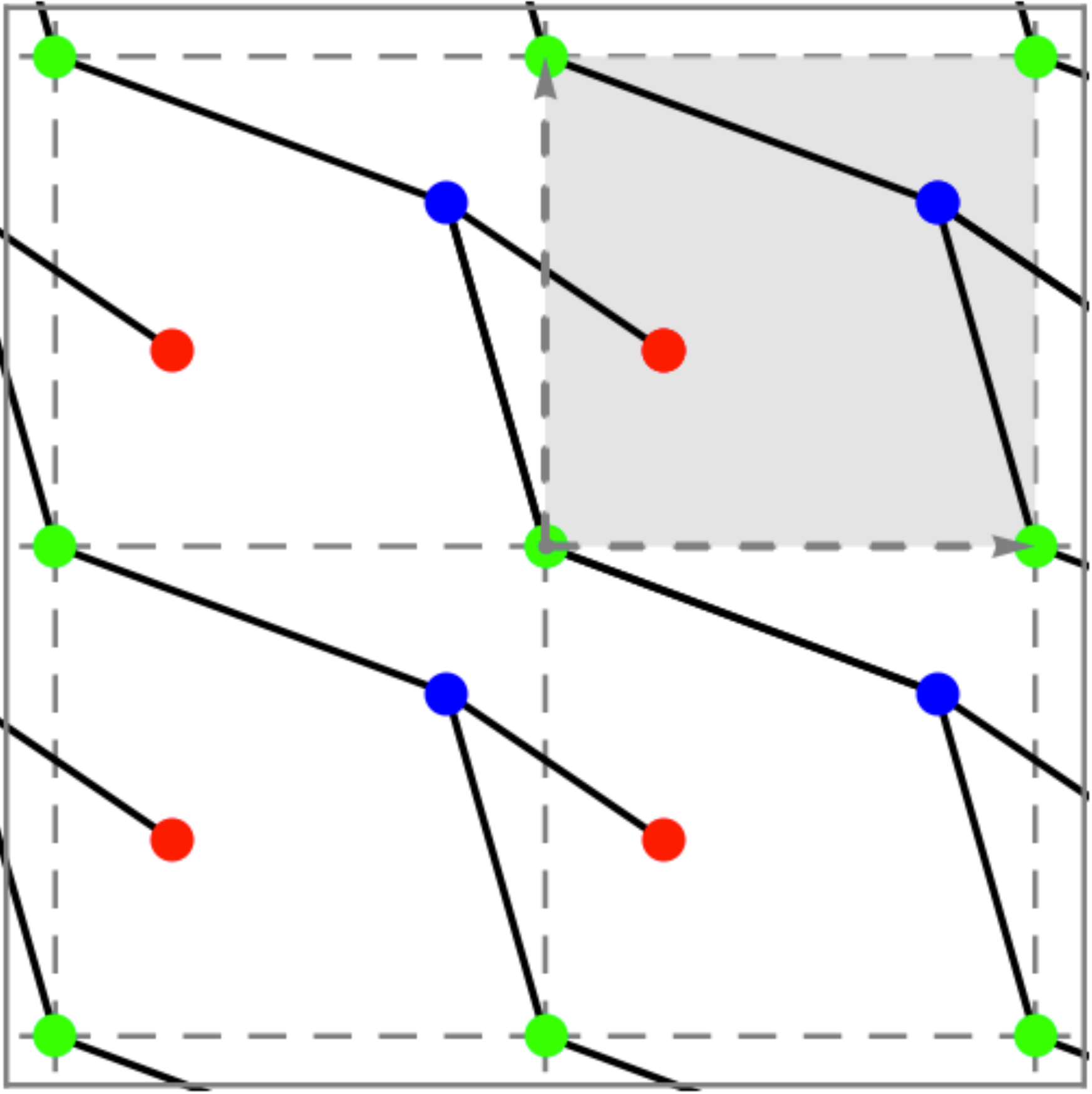}}\\
 \vspace{4pt}
 {\includegraphics[width=0.23\textwidth]{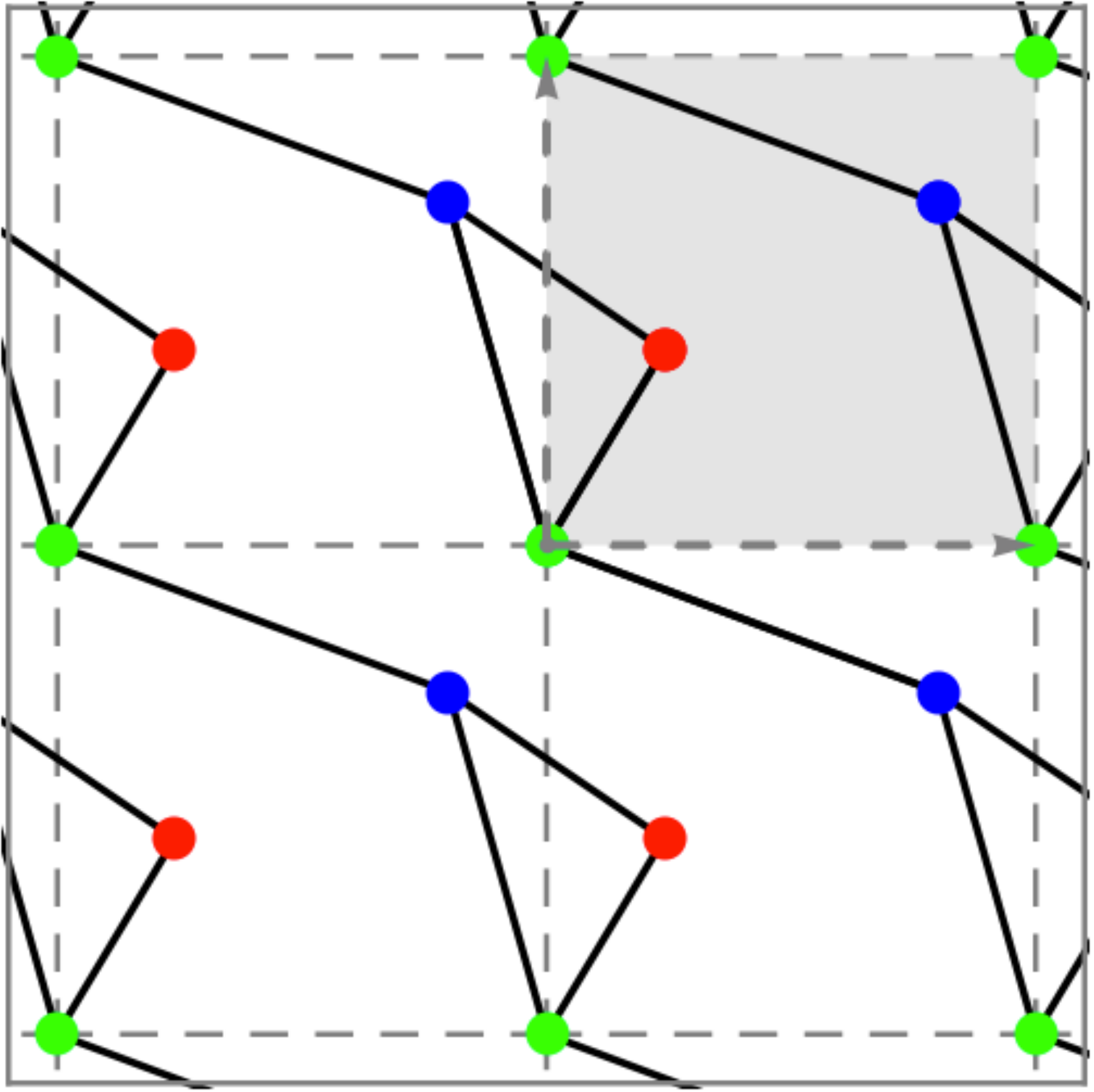}}\hspace{3pt}
 {\includegraphics[width=0.23\textwidth]{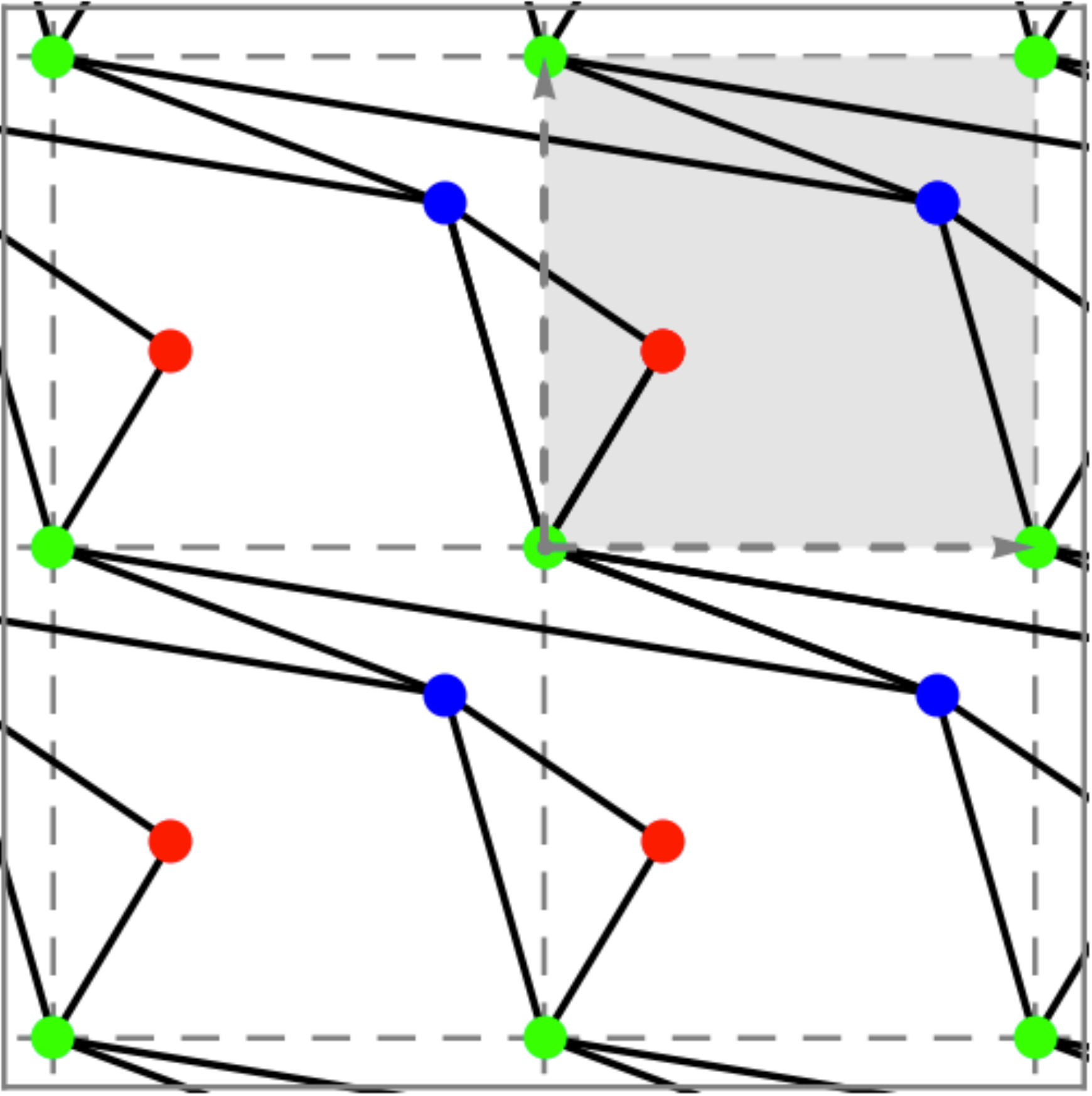}}\hspace{3pt}
 {\includegraphics[width=0.23\textwidth]{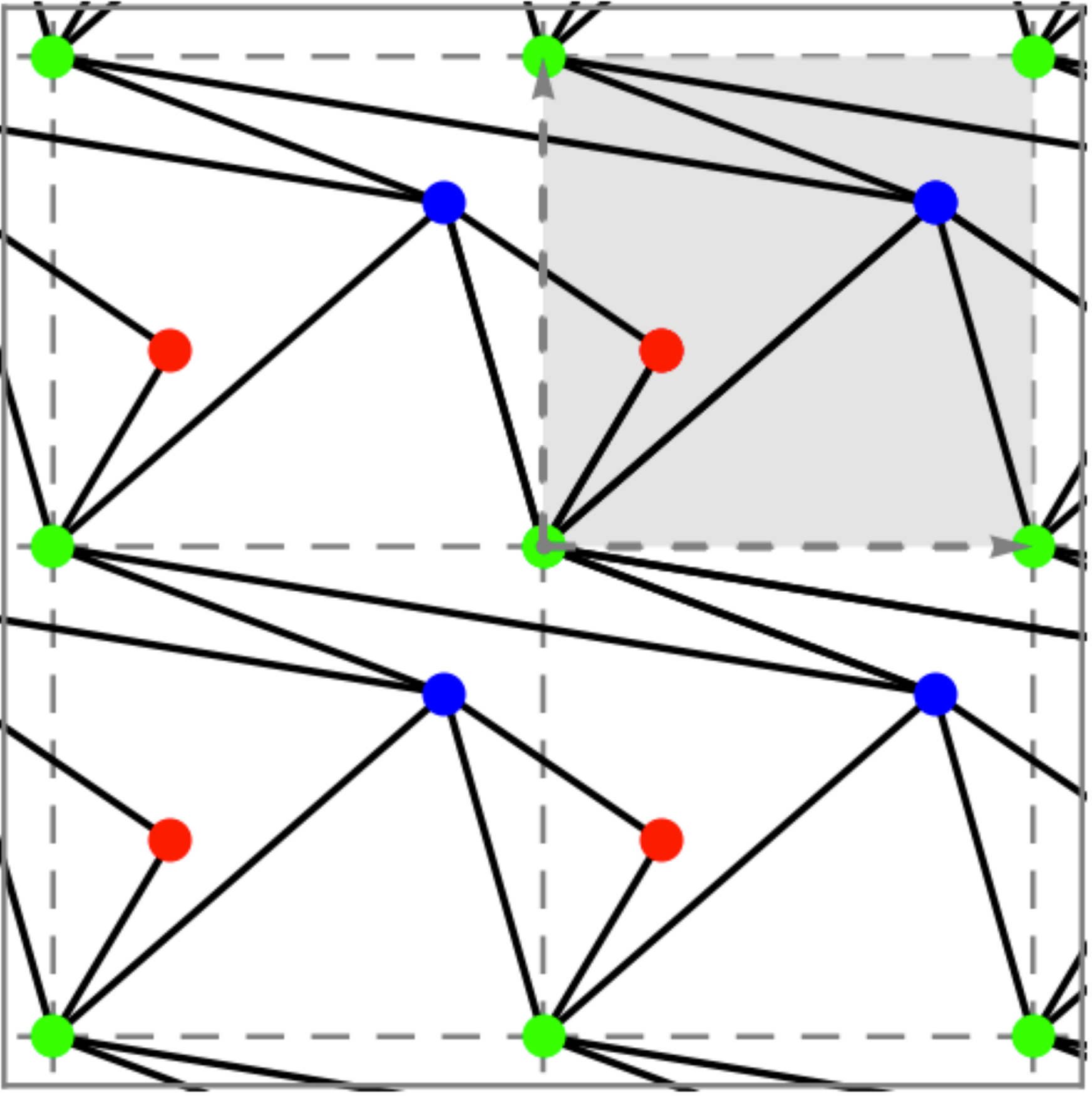}}\hspace{3pt}
 {\includegraphics[width=0.23\textwidth]{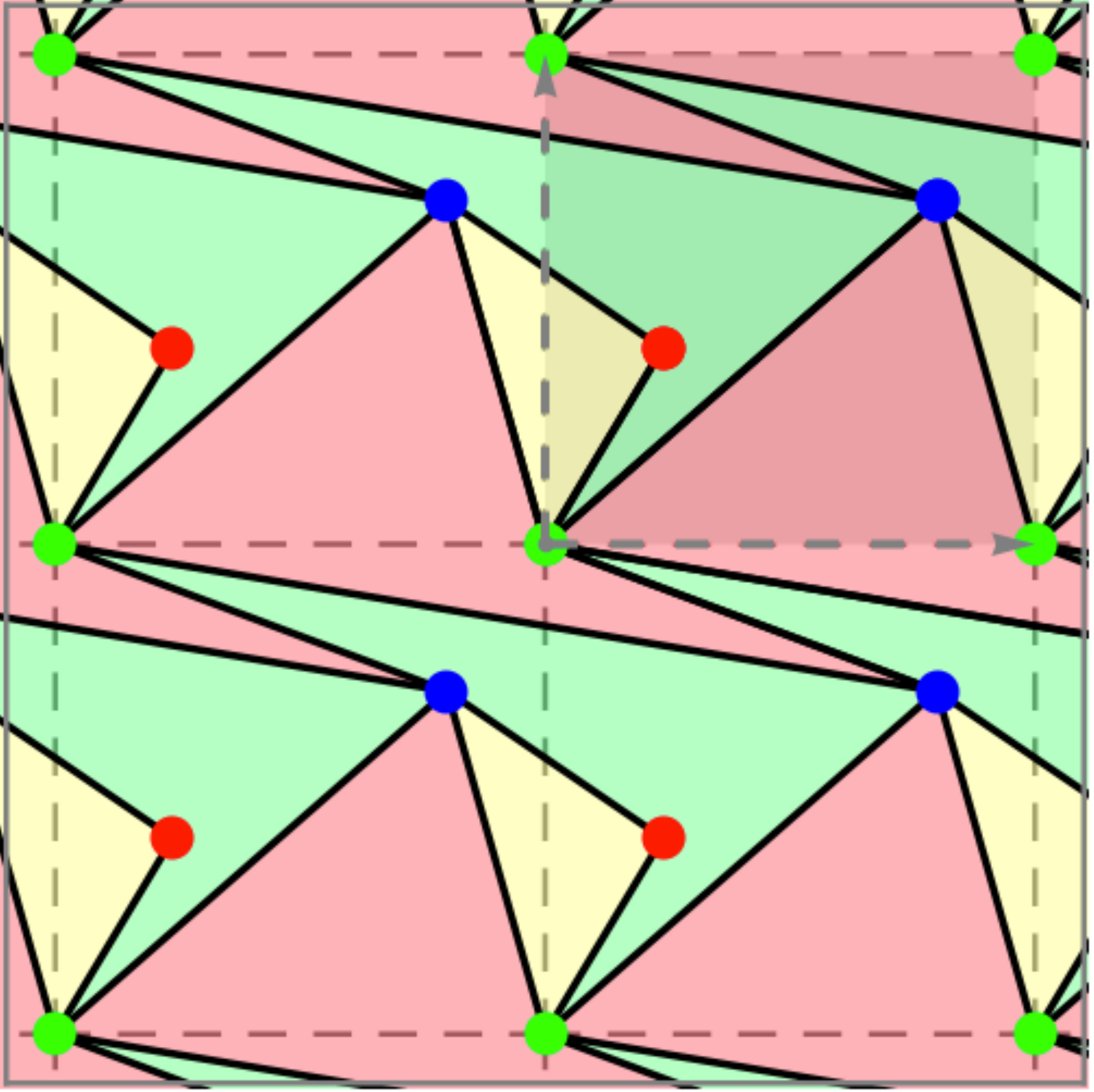}}
 \caption{  Generating a planar periodic pseudo-triangulation on $3$ fixed vertex orbits (colored red, blue and green). At each step, a random edge is inserted (and repeated by periodicity), subject only to maintaining the pointedness and non-crossing property of the framework. The faces of the final tiling are colored to emphasize their pseudo-triangular shape.}
 \label{FigPPT}
\end{figure}

\noindent
{\bf Infinitely many planar auxetic designs.}
We have shown above the far-reaching role of the notions of pseudo-triangulation and expansive behavior for planar periodic frameworks and auxetic investigations. We describe now a simple procedure for generating  periodic pseudo-triangulations {\em ad libitum}. Since expansive implies auxetic, we obtain {\em an infinite collection of auxetic designs.}

\medskip
\noindent
The generating procedure starts with an arbitrarily chosen lattice of periods and an arbitrarily chosen number $n$ of vertex orbits, as in Fig.~\ref{FigPPT}(a). Thus, the initial stage has $n$ orbits of vertices placed in the plane. From here, every new step will consist in inserting one edge orbit, subject to just two requirements: (i) to maintain pointedness at every vertex and (ii) not to produce edge crossings. It is proven in \cite{BS6} that exactly $m=2n$ steps are possible, and the end result is always a periodic pseudo-triangulation. Fig.~\ref{FigPPT} illustrates the procedure for the standard square lattice of periods and a placement of $n=3$ orbits of vertices.

\medskip 
\noindent
It may be observed that the class of periodic pseudo-triangulations is invariant under affine transformations and adoption of the standard square periodicity lattice at the start of the procedure is thereby warranted.

\medskip 
\noindent
Finally, we recall that, in general, auxetic behavior does not require pointedness. Hence many other types of auxetic designs are possible. For example, the rotating triangles of \cite{GAE} can be verified to satisfy our definition of auxeticity, but they are not, obviously, expansive. The question of verifying whether a given design is auxetic will be addressed in Section \ref{sec:scenario}.


 \section{The Minkowski space of symmetric $2\times 2$  matrices: auxetic trajectories as causal lines}
\label{sec:Minkowski}

In the planar case $d=2$, the analogy between auxetic trajectories for periodic frameworks and causal trajectories in a Minkowski space-time of dimension three turns into literal matching. Indeed, the three-dimensional vector space of symmetric $2\times 2$ matrices:
$$ A=A^t=(a_{ij})_{1\leq i,j\leq 2}, \ \ a_{12}=a_{21}  $$
\noindent
has a natural quadratic form: 
\begin{equation}\label{eq:det}
A\equiv (a_{11},a_{12},a_{22}) \mapsto det(A)=a_{11}a_{22}-a_{12}^2
\end{equation}
\noindent
This quadratic form has signature $(2,1)$, since: 
$$ a_{11}a_{22}-a_{12}^2=\frac{1}{2}(a_{11}+a_{22})^2-\frac{1}{2}(a_{11}-a_{22})^2-a_{12}^2 $$

\begin{wrapfigure}{l}{0.3\textwidth}
\vspace{-10pt}
\centering
\includegraphics[width=0.14\textwidth]{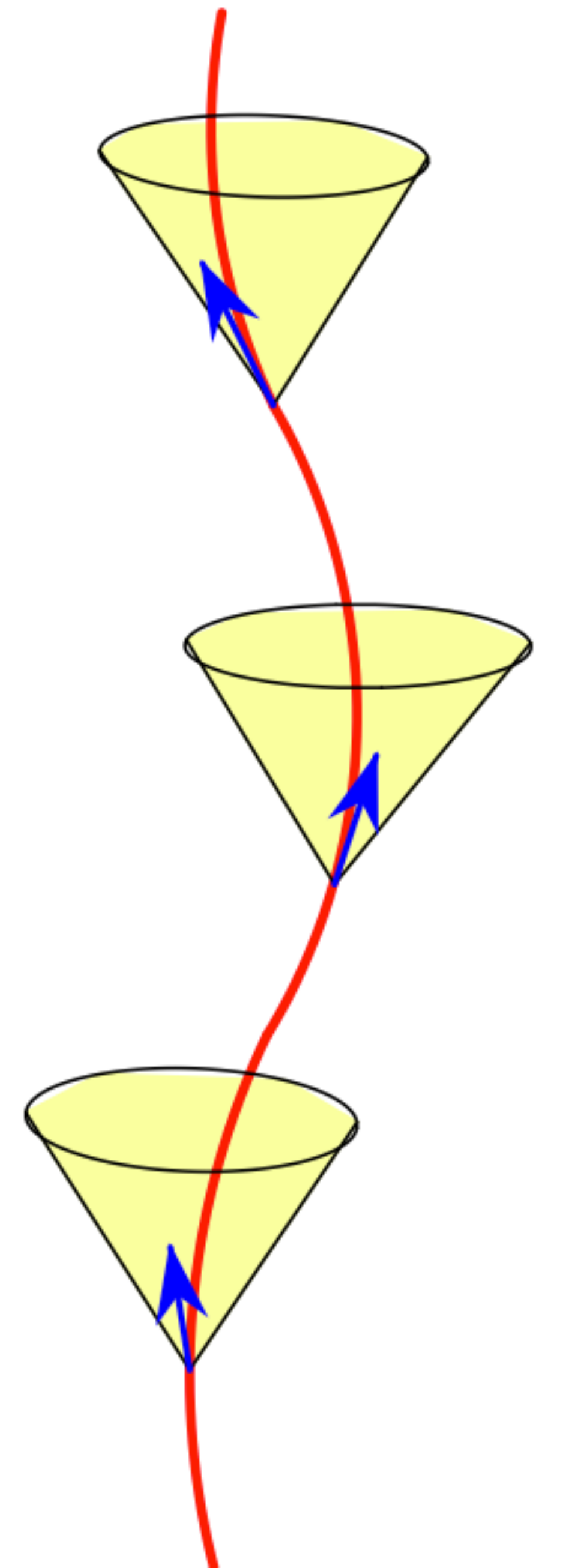}
 \caption{ Auxtic trajectories as causal lines.}
 \label{FigCausalLine}
\vspace{-16pt}
\end{wrapfigure}

\noindent
and defines a structure of Minkowski space-time, with space-like vectors of negative squared norm $det(A) < 0$ (Fig.~\ref{FigCausalLine}). The  `light cone' (at the origin) is defined by $det(A)=0$ and `future oriented' time-like vectors inside this cone are those with $a_{11}>0$ and $det(A)>0$ i.e. precisely the positive definite symmetric matrices. Thus, by Theorem~\ref{thm:tangent},  auxetic trajectories correspond precisely with causal trajectories (traced by the Gram matrices of a chosen pair of periodicity generators). One example is depicted in Fig.~\ref{FigCausalLine}.

\medskip
\noindent
We use this setting to gain further insight into  auxetic capabilities for hexagonal honeycomb frameworks. For simplicity, we assume hexagons  with equal edges and denote their common squared length by $s$. A concise manner to specify a periodic framework of this type is shown in Fig.~\ref{FigHexagonal}. It is enough to give the placement of a vertex representing one vertex orbit and show the three bars connecting it to three vertices in the other vertex orbit; the periodicity lattice is generated by the edge vectors of this triangle of vertices. With bars of the same length, the framework description amounts to a triangle of period vectors with the center of the circumscribed circle connected by bars to the vertices.
 
\begin{figure}[h]
\centering
 {\includegraphics[width=0.3\textwidth]{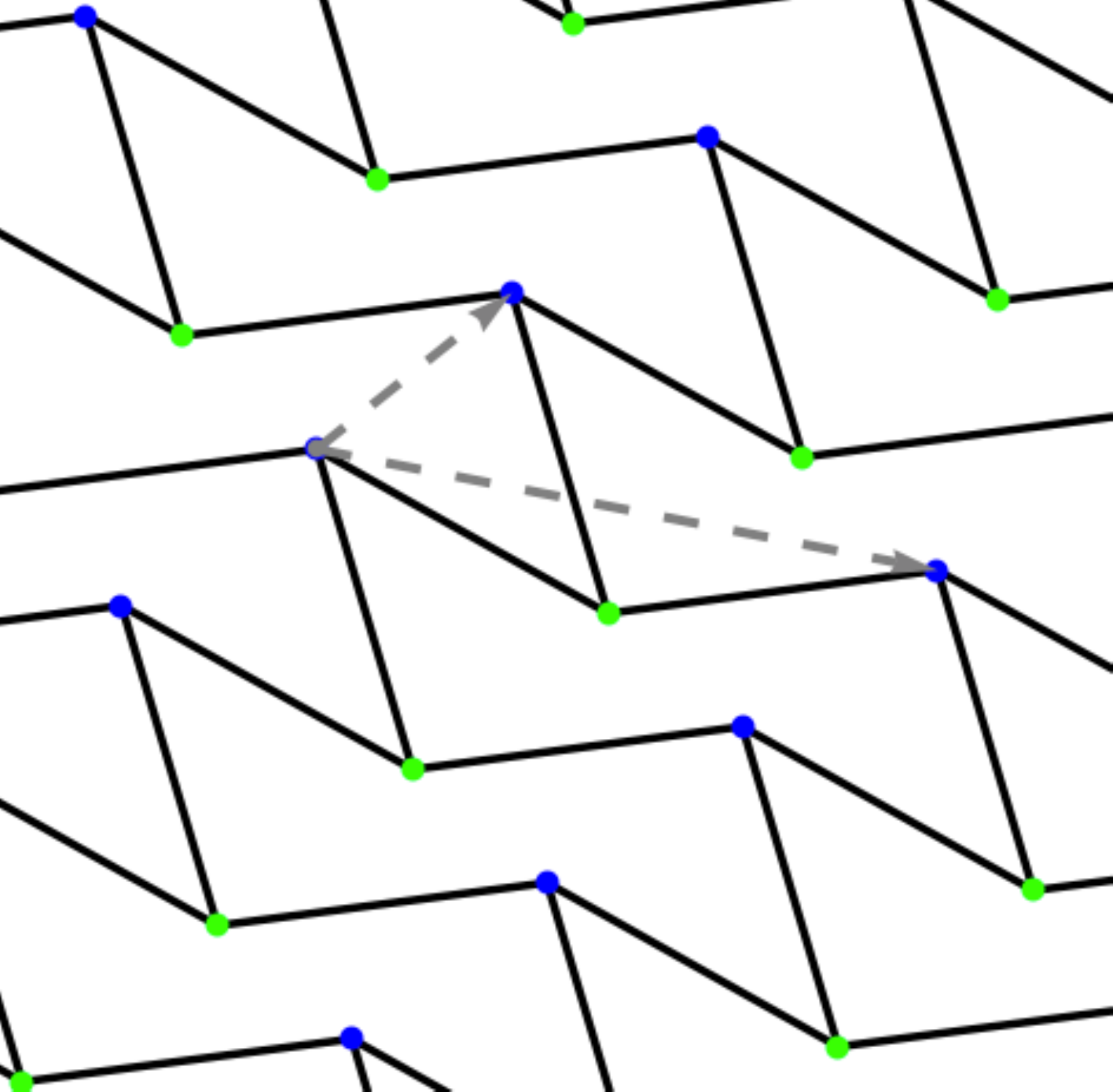}}
 \hspace{20pt}
 {\includegraphics[width=0.3\textwidth]{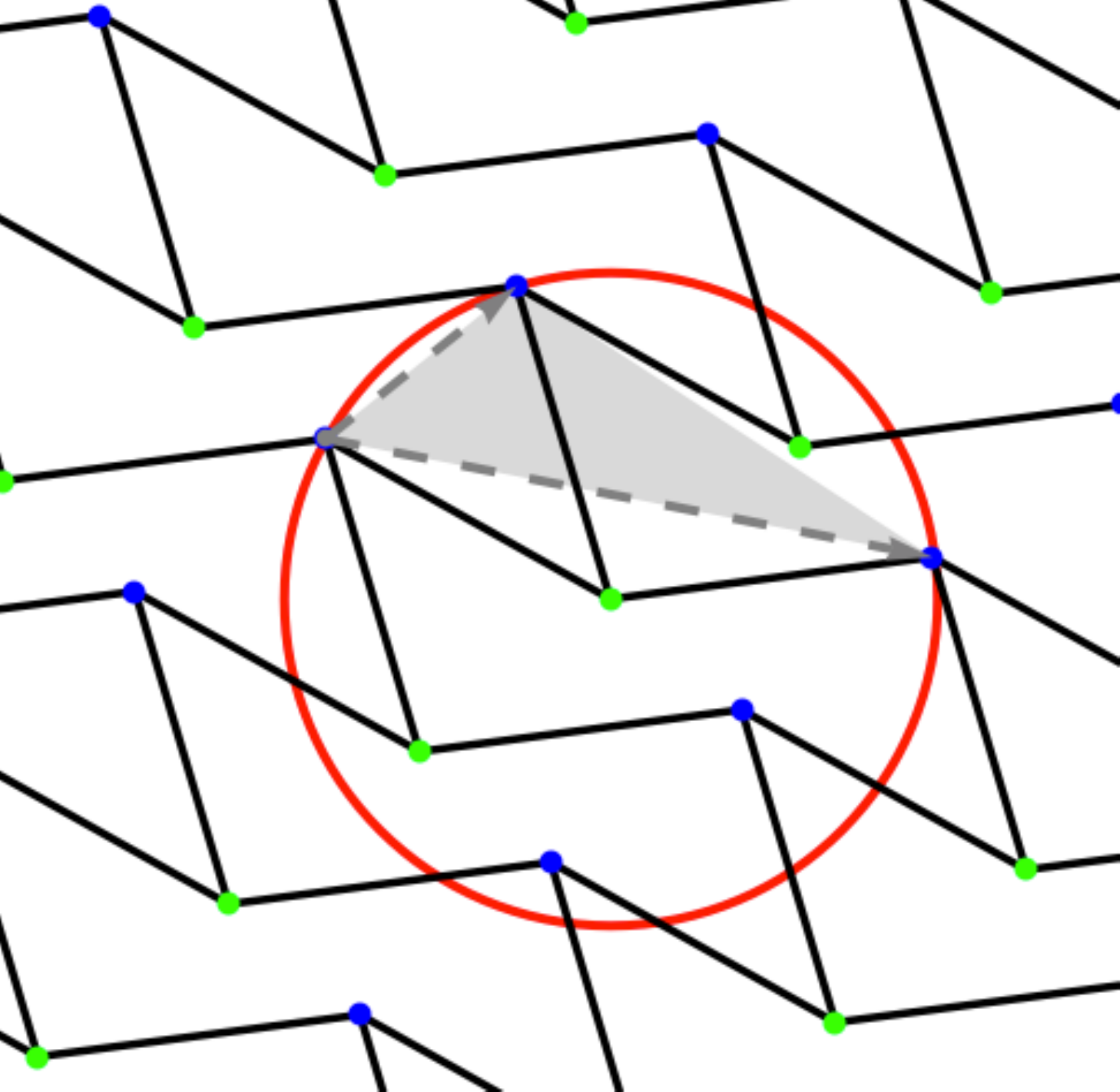}}
 {\includegraphics[width=0.3\textwidth]{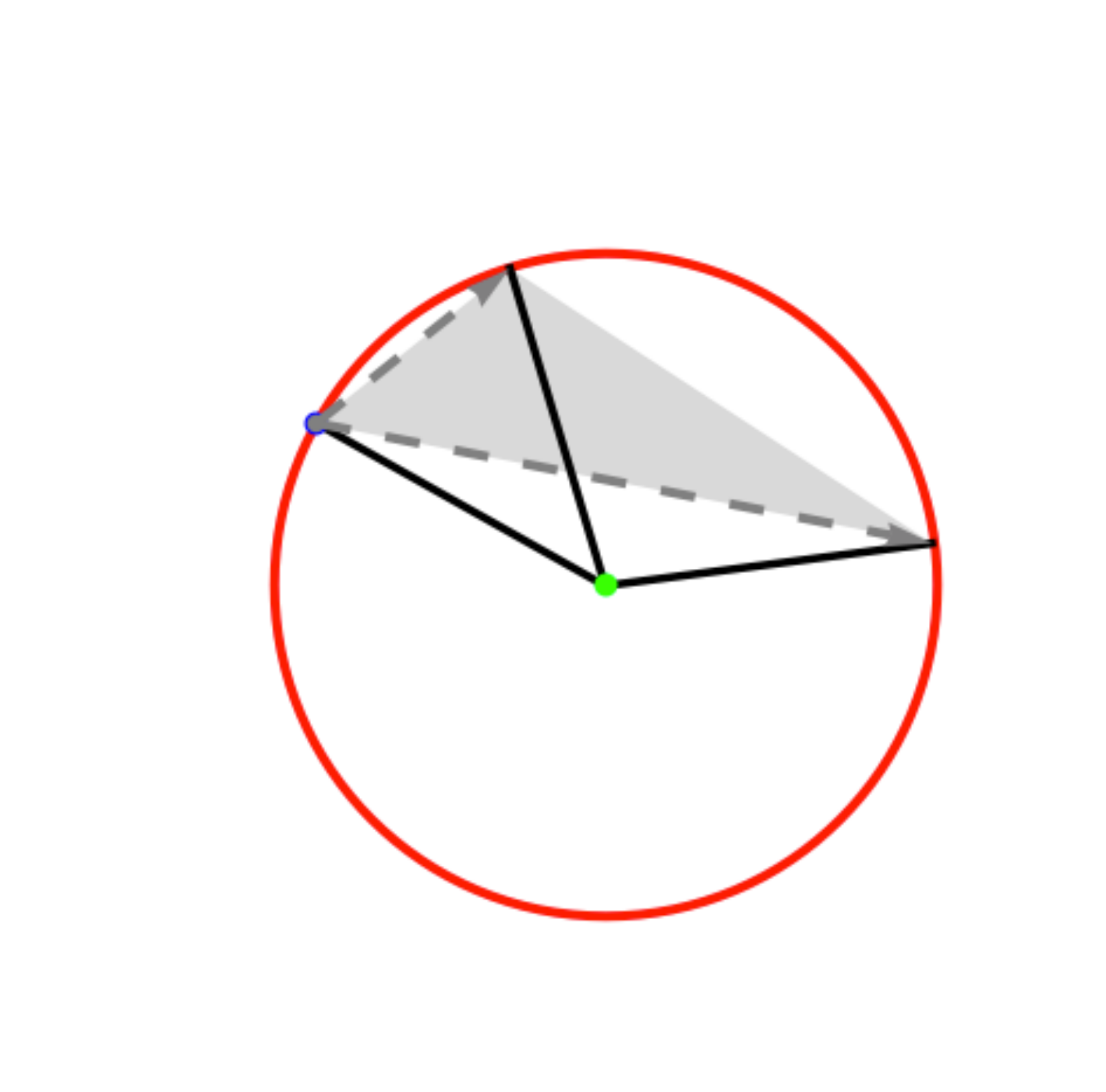}}
 \caption{The hexagonal honeycomb has auxetic deformations precisely when pointed at all vertices.}
 \label{FigHexagonal}
\end{figure}

\medskip
\noindent
The generators of the periodicity lattice are, as in Fig.~\ref{FigHexagonal}, $\lambda_1$ and $\lambda_2$,  with the corresponding Gram matrix entries denoted by $a_{ij}=\langle \lambda_i,\lambda_j \rangle$. The center of the circumscribed circle is at $\sigma$, with $\langle \sigma,\sigma \rangle=s$ fixed. A classical formula for the radius leads directly to the equation of the deformation space of the framework as a surface in the three-dimensional space with coordinates $(a_{11},a_{12},a_{22})$, namely:
\begin{equation}\label{eq:circum}
\frac{a_{11}a_{22}(a_{11}+a_{22}-2a_{12})}{4(a_{11}a_{22}-a_{12}^2)}=s
\end{equation}
%
\noindent
In order to determine the tangent plane at a given point of this surface, we rewrite
the defining equation as the (affine) cubic:
\begin{equation}\label{eq:cubic}
f(a_{11},a_{12},a_{22})=a_{11}a_{22}(a_{11}+a_{22}-2a_{12})-4s(a_{11}a_{22}-a_{12}^2)=0
\end{equation}
\noindent
and compute the components of the gradient $\bigtriangledown f$ of $f$:

\smallskip
\noindent
\hspace{14pt} $ f_{11}=\frac{\partial f}{\partial a_{11}}=a_{22}(a_{22}+2a_{11}-2a_{12}-4s) $

\smallskip
\noindent
\hspace{14pt} $ f_{12}=\frac{\partial f}{\partial a_{12}}=8sa_{12}-2a_{11}a_{22} $

\smallskip
\noindent
\hspace{14pt} $ f_{22}=\frac{\partial f}{\partial a_{22}}=a_{11}(a_{11}+2a_{22}-2a_{12}-4s) $

\medskip
\noindent
A tangent vector $\alpha=(\alpha_{11},\alpha_{12},\alpha_{22})$ must satisfy 
$\langle \alpha, \bigtriangledown f \rangle=0$, hence

\begin{equation}\label{eq:Tplane}
 \alpha_{12}= - \frac{\alpha_{11}f_{11}+\alpha_{22}f_{22}}{f_{12}} 
\end{equation}

\medskip
We have reached the point where we can address the question: {\em are there (non-trivial infinitesimal) auxetic deformations for the given framework?}

\medskip 
\noindent
The answer is obtained from the following considerations. Since an infinitesimal auxetic deformation must lie in the light cone $\alpha_{11}\alpha_{22}-\alpha_{12}^2\geq 0$, our tangent plane (\ref{eq:Tplane}) must intersect non-trivially the boundary cone $\alpha_{11}\alpha_{22}-\alpha_{12}^2= 0$. This means that we must have non-zero real solutions for the (homogeneous) quadratic equation:
\begin{equation}\label{eq:eq}
 \alpha_{11}\alpha_{22}= \left( \frac{\alpha_{11}f_{11}+\alpha_{22}f_{22}}{f_{12}} \right)^2
\end{equation}
\noindent
The discriminant takes the form:
\begin{equation}\label{eq:discriminant}
\Delta=(2f_{11}f_{22}-f_{12}^2)^2-4f_{11}^2f_{22}^2
\end{equation}
%
\noindent
For brevity, we omit the algebraic details of computing (and factoring) the discriminant of  equation (\ref{eq:eq}) as an expression in $(a_{11},a_{12},a_{22})$. Eventually, the condition $\Delta \geq 0$ has a simple and immediate geometrical interpretation: {\em  the triangle of periods must have an obtuse angle}. This means that the center of the circumscribed circle is not in the interior of the triangle and thus the three bars at the center are {\em pointed}. We obtain the following result.

\begin{prop}
\label{prop:hexagonal}
A necessary and sufficient condition for a planar periodic hexagonal framework (with equal edges) to allow non-trivial infinitesimal auxetic deformations is that of pointedness at all vertices.
\end{prop}

\medskip\noindent
An alternative derivation of this result appears in \cite{BS5}. For the type of framework under consideration, this clarifies the distinction between auxetic and expansive capabilities: for a non-trivial infinitesimal auxetic cone we need pointedness at all vertices, while for a non-trivial infinitesimal expansive cone we need pointedness and non-crossing.

\section{The general scenario in arbitrary dimension}
\label{sec:scenario}

In this section we return to arbitrary dimension $d$ and show that our theory can address the general problem of detecting auxetic trajectories in the local deformation space of a given $d$-periodic framework.
 
\medskip 
\noindent
The key elements of our purely mathematical approach are the following: (i) the notion of $d$-periodic framework  ${\cal F}=(G,\Gamma,p,\pi)$ and its {\em deformation space} $D({\cal F})$, (ii) when the deformation space has positive (i.e. non-zero) dimension, the notion of {\em auxetic trajectory}, defined for curves through the given framework in its deformation space, (iii) the vector space $Sym(d)$ of $d\times d$ symmetric matrices, with its {\em positive semidefinite cone} $\bar{\Omega}(d)\subset Sym(d)$ made of symmetric matrices with non-negative eigenvalues.

\medskip \noindent
These key elements are connected in a natural way. There is a map from the deformation space $D({\cal F})$ of the framework to the  vector space $Sym(d)$ obtained (after a choice of $d$ generators for the periodicity group $\Gamma$) by taking the Gram matrix of the corresponding generators of the periodicity lattice (which varies with the deformation).
\begin{equation}\label{eq:gMap}
 g:  D({\cal F}) \rightarrow Sym(d)  
\end{equation} 

\noindent
Actually the image is contained in the interior of the positive semidefinite cone. Thus, (parametrized) curves through ${\cal F}\in D({\cal F})$ are mapped to (parametrized) curves in $Sym(d)$ and we can distinguish velocity vectors which belong (as free vectors) to the positive semidefinite cone $\bar{\Omega}(d)$ as the infinitesimal mark of an auxetic tendency. In other words, {\em a one-parameter deformation is auxetic when all its velocity vectors belong to the positive semidefinite cone.}

\medskip \noindent
Thus, for our framework ${\cal F}$, we can look at the {\em tangent space} at ${\cal F}$ to the deformation space and identify as {\em infinitesimal auxetic deformations} those which  map (by the differential $T(g)$ of $g$) to vectors in the positive semidefinite cone. Hence, at every point of $D({\cal F})$, we have an infinitesimal {\em auxetic cone} in the corresponding tangent space and a curve in $D({\cal F})$ is an auxetic deformation when all its tangent vectors are in the respective auxetic cones (along the curve).

\begin{wrapfigure}{l}{0.5\textwidth}
\vspace{-18pt}
\centering
\includegraphics[width=0.5\textwidth]{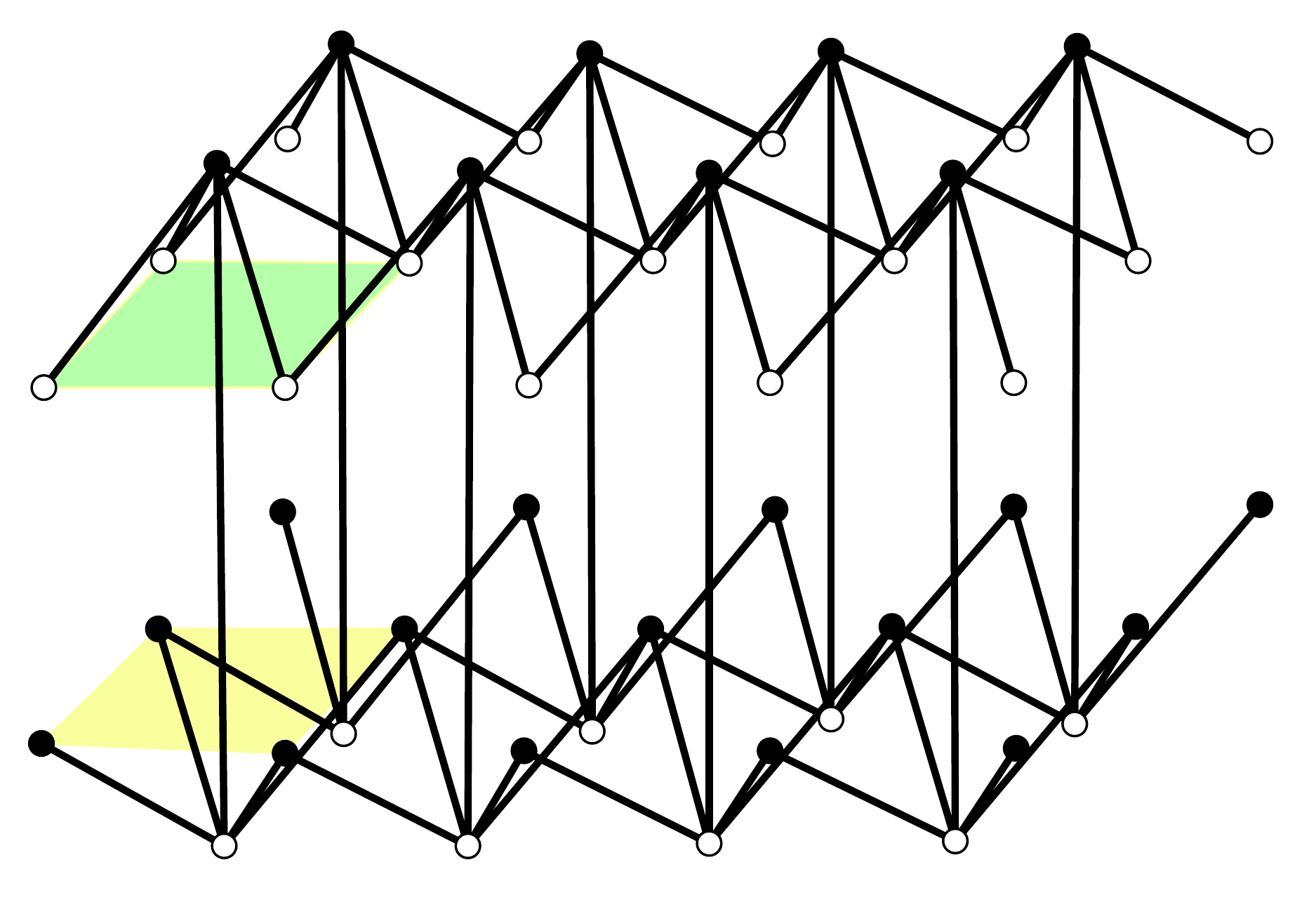}
\vspace{-15pt}
 \caption{ A 3D periodic framework with expansive and auxetic capabilities. }
 \label{FigExpansive4dof}
\vspace{-18pt}
\end{wrapfigure}

\medskip \noindent
{\bf Remarks.}\ Since all tangent spaces for points in $Sym(d)$ are identified with $Sym(d)$, the tangent map for $g$ in  (\ref{eq:gMap}) gives, for the framework ${\cal F}$:
\begin{equation}\label{eq:TgMap}
 T_{\cal F} (g):  T_{\cal F} D({\cal F}) \rightarrow Sym(d)  
\end{equation} 
\noindent
In $Sym(d)$ we have to look at the intersection of the image of $T_{\cal F} (g)$ (a vector subspace) with the positive semidefinite cone $\bar{\Omega}(d)$. This is known in polynomial optimization and semidefinite programming as a {\em spectrahedral cone} \cite{Vin}. In this language, infinitesimal auxetic cones \ are \ linear \ preimages \ of spectrahedral cones. 

\medskip \noindent
From this scenario, it becomes evident that, when a framework has sufficiently many degrees of freedom, it is fairly likely to possess some auxetic capabilities, that is, certain auxetic one-parameter deformations. In particular, if ${\cal F}\in D({\cal F})$ is a smooth point and  $T_{\cal F}(g)$ in (\ref{eq:TgMap}) is onto, then auxetic trajectories are guaranteed locally.

\medskip
\noindent
{\bf Deciding infinitesimal auxeticity.}
In the manner described above, the problem of deciding if a given periodic framework has a non-trivial auxetic infinitesimal deformation turns into a feasibility problem in semidefinite programming and can be addressed computationally using well-studied algorithms and several available, efficient implementations \cite{BPT,PK}. In particular, since the question is posed for a fixed dimension $d$, the algorithm of \cite{PK} runs in polynomial time. We obtain:
\begin{corollary}
\label{cor:decideAuxeticity}
The problem of deciding if a $d$-periodic framework allows an infinitesimal auxetic deformation can be decided in polynomial time using semidefinite programming.
\end{corollary}

\medskip
\noindent
{\bf A case study in dimension three.}
We apply the above scenario to the three dimensional periodic framework illustrated in Fig.~\ref{FigExpansive4dof}. The framework has $n=2$ vertex orbits and $m=5$ edge orbits, resulting in 4 degrees of freedom.  The essential elements for constructing this periodic framework and its deformations are shown in Fig.~\ref{FigBasics4dof}(a). The periodicity lattice is generated by the three vectors:
\begin{equation}\label{eq:gen}
\overrightarrow{OP}= \left( \begin{array}{c} 2\alpha_1\\ 
0 \\
0  \end{array} \right),\ 
\overrightarrow{OR}= \left( \begin{array}{c}  0\\
2\alpha_2\\ 
 
0  \end{array} \right),\ \ 
\overrightarrow{OS}=\beta= \left( \begin{array}{c} \beta_1\\ 
\beta_2 \\
\beta_3  \end{array} \right). 
\end{equation}
\noindent
The four edges $AO,AP,AQ,AR$ are taken of equal length 1, with the fifth edge $AS$ of length
$r=3/\sqrt{5}$. Since the (period) parallelogram $OPQR$ remains under deformations an inscribed parallelogram, it remains rectangular. We keep $O$ fixed as the origin. Thus:
\begin{equation}\label{eq:param}
\overrightarrow{OA}=\alpha=\left( \begin{array}{c} \alpha_1\\ 
\alpha_2 \\
\alpha_3  \end{array} \right),\ \mbox{with}\ \ \langle \alpha,\alpha \rangle=1
\end{equation}
\noindent
The framework is completely described by the six parameters in $\alpha, \beta$, constrained by
two relations, namely:
\vspace{-8pt}
\begin{equation}\label{eq:rel}
 \langle \alpha,\alpha \rangle=1 \ \ \mbox{and} \ \ 
\langle \beta -\alpha, \beta -\alpha \rangle=r^2=\frac{9}{5}
\end{equation}

\vspace{-18pt}
\begin{figure}[h]
\centering
 \subfloat[]{\includegraphics[width=0.4\textwidth]{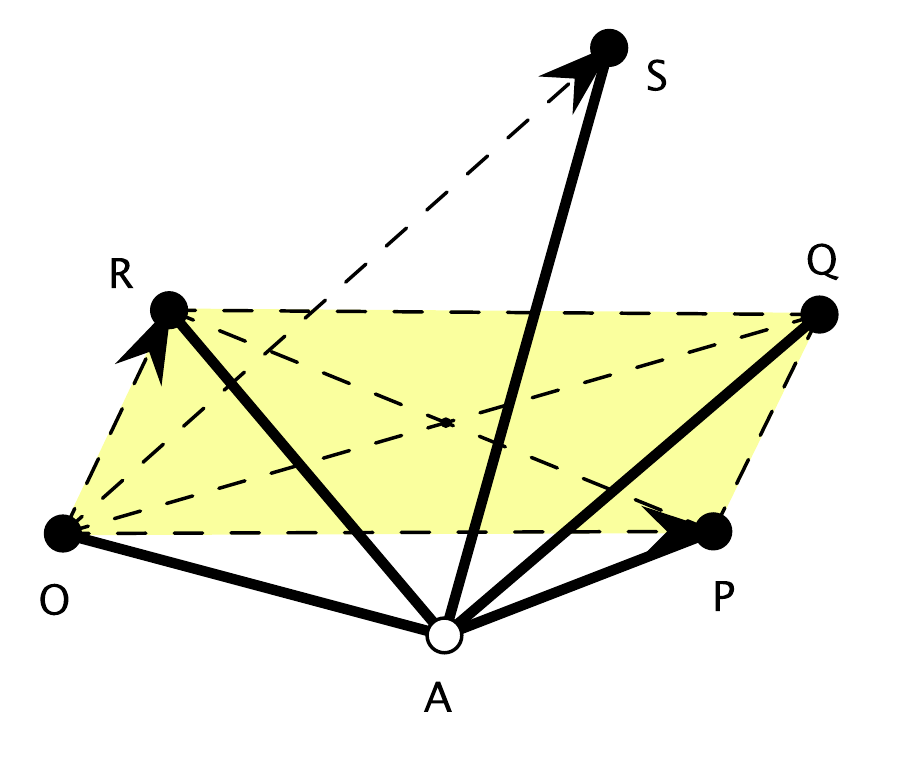}}
 \subfloat[]{\includegraphics[width=0.6\textwidth]{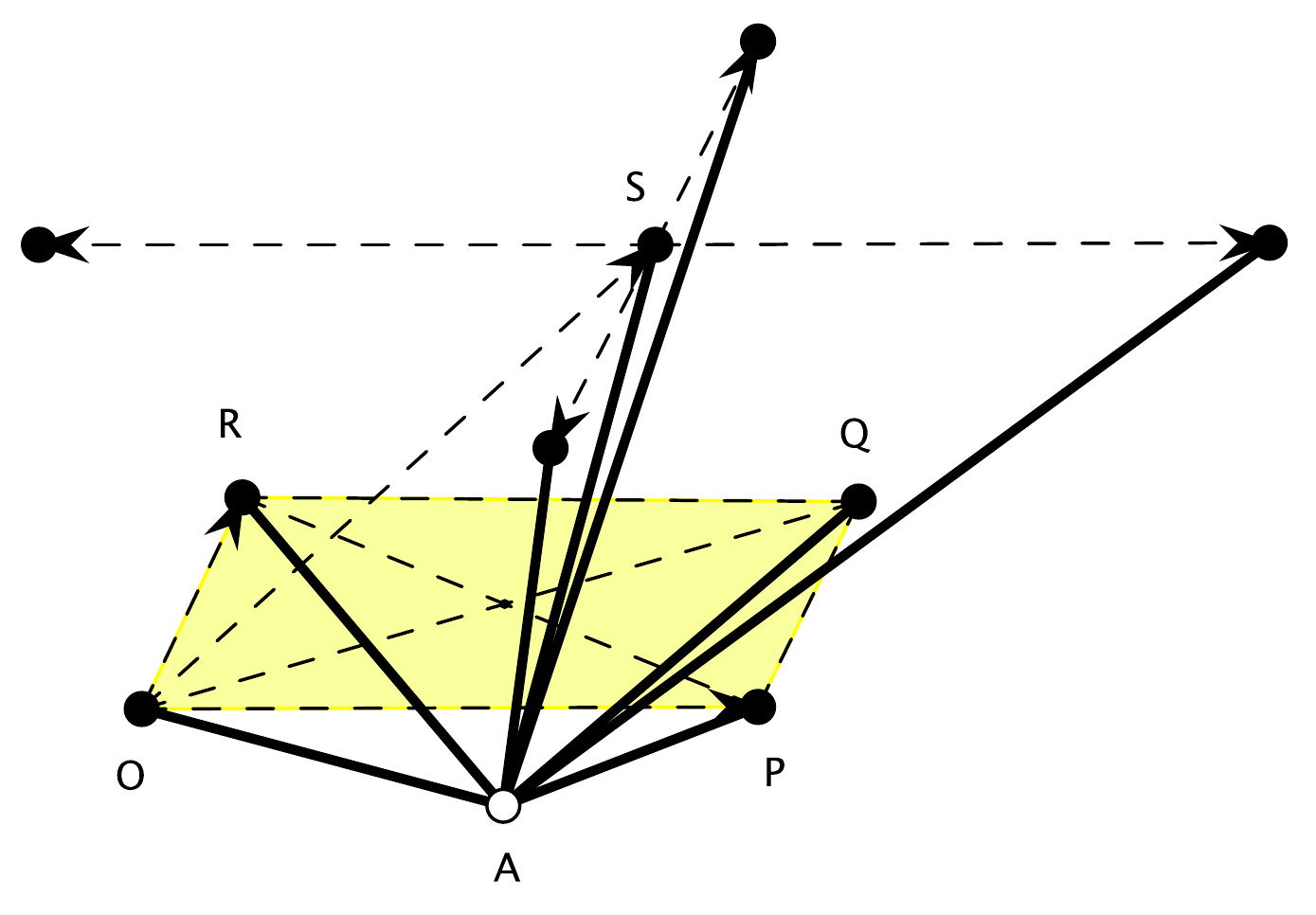}}
 \caption{ (a) The essential information for the analysis of the 3D framework from Fig.\ref{FigExpansive4dof}. (b) Inserting three new edge orbits.}
 \label{FigBasics4dof}
 \label{FigBasics1dof}
\end{figure}

\noindent
We adopt as initial position the configuration where $OPQRS$ is the upper half of a regular octahedron of squared edge $8/5$, that is: $\alpha_1=\alpha_2=\beta_1=\beta_2=\sqrt{2/5}$ and $\alpha_3=-1/\sqrt{5}$, $\beta_3=2/\sqrt{5}$. The Gram matrix map (\ref{eq:gMap}) is described in parameters $\alpha,\beta$ by:
\begin{equation}\label{eq:Gram1}
\omega  =
\left( \begin{array}{lll} 
	4\alpha_1^2 & 0 & 2\alpha_1\beta_1 \\
    0 & 4\alpha_2^2 & 2\alpha_2\beta_2 \\
    2\alpha_1\beta_1 & 2\alpha_2\beta_2 & \beta_1^2+\beta_2^2+\beta_3^2 
	\end{array}  
\right)=
\left( \begin{array}{lll}  
	a_{11} & 0 & a_{13}\\
    0 & a_{22} & a_{23} \\
    a_{13} & a_{23} & a_{33} 
\end{array}  \right)
\end{equation}
\noindent
The two relations in (\ref{eq:rel}) imply the already noticed orthogonality $a_{12}=0$ and the quartic equation:
\begin{equation}\label{eq:quartic}
f( a ) = a_{11}a_{22}(a_{33}-a_{13}-a_{23}+1-r^2)^2 -\Delta (4-a_{11}-a_{22}) =0
\end{equation}
\noindent
where $\Delta=det( \omega )=a_{11}a_{22}a_{33}-a_{22}a_{13}^2-a_{11}a_{23}^2$. For a  description of the  {\bf local deformation space} of our framework we may use directly the quartic hypersurface (\ref{eq:quartic}) in the space of five variables $a=(a_{11},a_{22},a_{33},a_{13},a_{23})$, in a neighborhood of the initial position $a(0)=(8/5,8/5,8/5,4/5,4/5)$.

\smallskip
\noindent
The components of the gradient $\bigtriangledown(f)(a)$ are as follows:

\smallskip \noindent
$ f_{11}=\frac{\partial f}{\partial a_{11}}=a_{22}(a_{33}-a_{13}-a_{23}+1-r^2)^2 
-(a_{22}a_{33}-a_{23}^2)(4-a_{11}-a_{22}) + \Delta $

\smallskip \noindent
$ f_{22}=\frac{\partial f}{\partial a_{22}}=a_{11}(a_{33}-a_{13}-a_{23}+1-r^2)^2 
-(a_{11}a_{33}-a_{13}^2)(4-a_{11}-a_{22}) + \Delta $

\smallskip \noindent
$ f_{33}= \frac{\partial f}{\partial a_{33}}=2a_{11}a_{22}(a_{33}-a_{13}-a_{23}+1-r^2) 
-a_{11}a_{22}(4-a_{11}-a_{22})   $ 

\smallskip \noindent
$ f_{13}=\frac{\partial f}{\partial a_{13}} = -2a_{11}a_{22}(a_{33}-a_{13}-a_{23}+1-r^2) 
+2a_{22}a_{13}(4-a_{11}-a_{22}) $

\smallskip \noindent
$ f_{23}=\frac{\partial f}{\partial a_{23}} = -2a_{11}a_{22}(a_{33}-a_{13}-a_{23}+1-r^2) 
+2a_{11}a_{23}(4-a_{11}-a_{22}) $

\medskip
\noindent
and the {\bf infinitesimal auxetic cone} at $a$ will be given by the intersection of the orthogonal
space $\bigtriangledown(f)(a)^{\perp}$ with the positive semidefinite cone (restricted to our
five dimensional subspace $a_{12}=0$). The resulting type of {\em spectrahedral
cone} is familiar: it is the cone over the `terahedral part' of a nodal Cayley cubic surface, as
illustrated e.g. in \cite{ORSV}. 

\begin{wrapfigure}{l}{0.44\textwidth}
\vspace{-8pt}
\centering
{\includegraphics[width=0.44\textwidth]{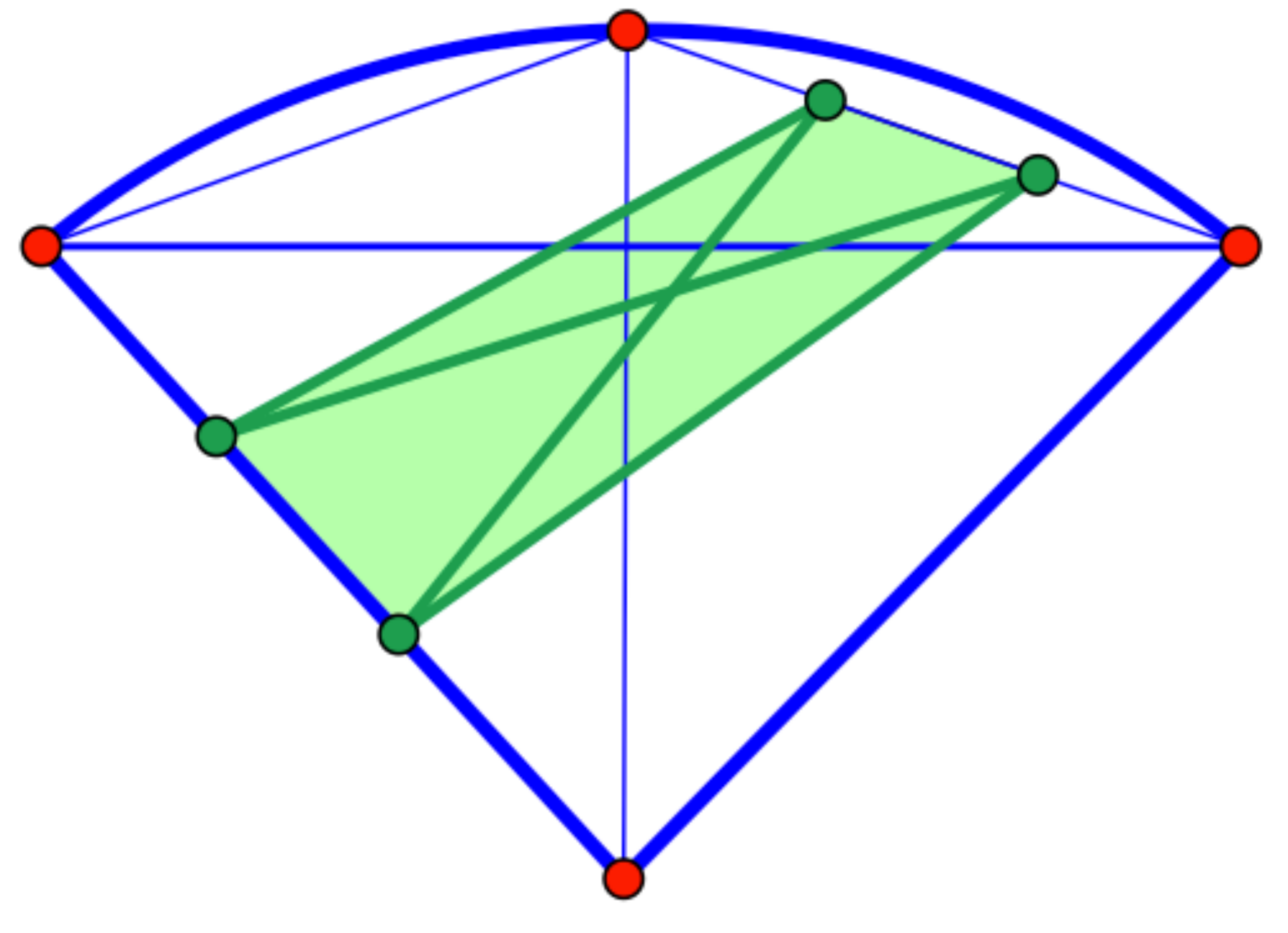}}
\vspace{-12pt}
 \caption{The polyhedral expansive cone included in the spectrahedral auxetic cone, shown  in a three dimensional section. The spectrahedron is bounded by a Cayley cubic whose four nodes correspond to rank one symmetric $3\times 3$ matrices.  }
 \label{FigAEcones}
\vspace{-24pt}
\end{wrapfigure}

\medskip
\noindent
Since our four-parameter deformation family is particularly apt to illustrate this qualitative feature, we present the necessary details. In order to refer to elements in the {\em tangent
bundle} of our quartic fourfold, we shall use pairs $(a,v)$ where $a$ is a point of our local deformation space identified with a neighborhood of $a(0)=(8/5,8/5,8/5,4/5,4/5)$ in the
quartic $f(a)=0$ of $R^5$, and $v$ is a vector in the tangent space at that point, identified with
$\bigtriangledown(f)(a)^{\perp}\subset R^5$.

\bigskip
\noindent
The boundary of the positive semidefinite cone $\bar{\Omega(3)}\subset Sym(3)\equiv R^6$
is contained in the homogeneous cubic locus $\Delta=0$, made of symmetric $3\times 3$ matrices of rank at most two. Thus, projectively, we are in $P_5$ and the {\em rank one locus}
corresponds to the quadratic Veronese embedding $P_2 \rightarrow P_5$, with image of
degree four \cite{B,Gr}. 
When intersected with the hyperplane $v_{12}=0$, this Veronese surface gives
two conics: one in  $v_{11}=v_{13}=0$, with equation $v_{22}v_{33}=v_{23}^2$ and the other in
$v_{22}=v_{23}=0$, with equation $v_{11}v_{33}=v_{13}^2$. 

\medskip
\noindent
When further intersecting with $\bigtriangledown(f)(a)^{\perp}$, that is with
$\langle \bigtriangledown(f)(a),v \rangle =0$, we obtain {\em four projective
points}. Using the abbreviations $\delta_{13}=(f_{13}^2-4f_{11}f_{33})^{1/2}$ and
$\delta_{23}=(f_{23}^2-4f_{22}f_{33})^{1/2}$, they are:
\begin{equation}\label{eq:nodes}
\begin{array}{l}
 (-f_{33}(\delta_{13}-f_{13}):0:f_{11}(\delta_{13}+f_{13}):-2f_{11}f_{33}:0) \\
 (-f_{33}(\delta_{13}+f_{13}):0:f_{11}(\delta_{13}-f_{13}):2f_{11}f_{33}:0)  \\
(0:-f_{33}(\delta_{23}-f_{23}):f_{22}(\delta_{23}+f_{23}):0:-2f_{22}f_{33}) \\
(0:-f_{33}(\delta_{23}+f_{23}):f_{22}(\delta_{23}-f_{23}):0:2f_{22}f_{33})
\end{array}
\end{equation}

\noindent
The cubic equation $\Delta(v)=v_{11}v_{22}v_{33}-v_{22}v_{13}^2-v_{11}v_{23}^2=0$ restricts to a Cayley cubic surface with its four nodes at the above points. The actual boundary of the
auxetic cone (in our projective description) retains only the portion of this Cayley cubic surface which `wraps around' the tetrahedron of nodes and contains its edges.

\medskip
\noindent
We proceed now with the determination of the {\bf infinitesimal expansive cone} for our framework, inside the spectrahedral auxetic cone identified above. For this purpose, it will be useful to consider the parallel plane to $OPQR$ though $S$ and mark the four periodicity translates $S_{\pm 1}= S\pm  \overrightarrow{OP}$ and $S_{\pm 2}= S\pm \overrightarrow{OR}$. We obtain four {\em one degree of freedom} mechanisms by inserting in our original framework three new edge orbits corresponding to three new bars connecting $A$ to three of the marked translates, as illustrated in Fig.~\ref{FigBasics1dof}(b).

\medskip
\noindent
If we denote by $r^2_{\pm k}$ the squared distances from $A$ to $S_{\pm k}$, $k=1,2$, we obtain the relations:
\begin{equation}\label{eq:4bars}
\begin{array}{l}                                    
2a_{11}= r^2_{-1}+r^2_{+1}-2r^2 \\
2a_{22}=r^2_{-2}+r^2_{+2}-2r^2 \\
2a_{13}=r^2_{+1}-r^2 \\
2a_{23}=r^2_{+2}-r^2
\end{array}
\end{equation}

\noindent
Now, we fix three of the parameters $r^2_{\pm k}$, $k=1,2$ and vary the fourth. This leads
easily to the infinitesimal deformations corresponding to our four mechanisms. They are:
\begin{equation}\label{eq:extremals}
\begin{array}{l}      
(-f_{33}:0:f_{11}:0:0) \\
(-f_{33}:0:f_{11}+f_{13}:-f_{33}:0) \\
(0:-f_{33}:f_{22}:0:0) \\
(0:-f_{33}:f_{22}+f_{23}:0:-f_{33})
\end{array}
\end{equation}

\noindent
The expansive and auxetic capabilities of our framework can be summarized as follows.

\begin{theorem}\label{thm:cones}
In a sufficiently small neighborhood of the initial position $a(0)=(8/5,8/5,8/5,4/5,4/5)$, the infinitesimal expansive cone and infinitesimal auxetic cone, parametrized by $a\in \{ f(a)=0\}$, of the framework in Fig.~\ref{FigExpansive4dof} correspond, projectively, to the tetrahedron with vertices at (\ref{eq:extremals}), respectively the spectrahedron with nodes at (\ref{eq:nodes}).
\end{theorem}

\noindent
{\em Proof:}\ The structure of the infinitesimal auxetic cone was described above as a cone over the portion of a Cayley cubic `wrapping around' the tetrahedron of nodes. The fact that, for deformations close enough to the initial position, all four mechanisms obtained from the indicated insertions of three new edge orbits are actually expansive can be directly verified. We omit further computational details, but illustrate the relative position of the polyhedral expansive cone included in the spectrahedral auxetic cone in Fig.~\ref{FigAEcones}.

\medskip
\noindent
For the initial framework, the gradient direction is $(1:1:-4:4:4)$, the nodes are at:
$$ (0:2(\sqrt{2}-1):(\sqrt{2}+1)/2:0:1) , \ \  (0:2(\sqrt{2}+1):(\sqrt{2}-1)/2:0:-1), $$
$$  (2(\sqrt{2}-1):0:(\sqrt{2}+1)/2:1:0), \ \  (2(\sqrt{2}+1):0:(\sqrt{2}-1)/2:-1:0) $$
\noindent
and the extremal rays of the infinitesimal expansive cone are given by
$$ (4:0:1:0:0),\ (4:0:5:4:0),\ (0:4:1:0:0), (0:4:5:0:4).\ \ \ \  \qed $$

\medskip \noindent
{\bf Remark.}\ Formulae (\ref{eq:nodes}) and (\ref{eq:extremals}) obtained above may be seen as specific direction fields apt to be combined by various linear combinations with positive functions as coefficients. Locally, integral curves for such direction fields will be auxetic and expansive trajectories respectively.

\section{New three-dimensional auxetic designs}
\label{sec:new3Ddesigns}

In this section we show that the three-dimensional structure studied above belongs to a family of kindred auxetic designs. While the framework type used in our case study was noticed in the earlier literature \cite{Alm}, we obtain new designs based on features necessarily
present in periodic frameworks with expansive capabilities. In \cite{BS7} we have identified a version of pointedness which must be present for effective expansiveness in dimensions higher than two and the family of examples presented here will illustrate the design significance of our {\em expansive implies auxetic} emphasis.

\medskip
\noindent
We describe six types of three-dimensional periodic frameworks, all with just {\em two orbits of vertices}. As recalled in Section~\ref{sec:deformations}, this means that, modulo periodicity, we have exactly two equivalence classes of vertices, rendered as black and white in Figs.~\ref{FigDesigns12},~\ref{FigDesigns34} and \ref{FigDesigns56}. If $(G,\Gamma)$ stands for the abstract $3$-periodic graph under consideration, with $G=(V,E)$, the cardinality $n = |V/\Gamma|$ remains $n=2$, but the number of edge orbits $m=|E/\Gamma|$ varies from one example to another.

\begin{wrapfigure}{l}{0.55\textwidth}
\vspace{-4pt}
\centering
{\includegraphics[width=0.24\textwidth]{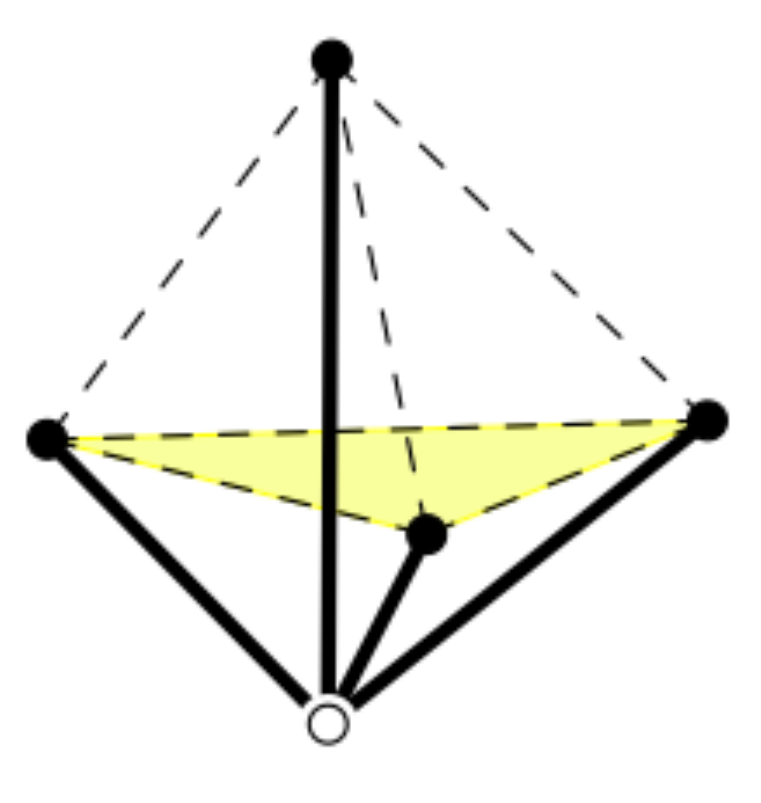}}
{\includegraphics[width=0.3\textwidth]{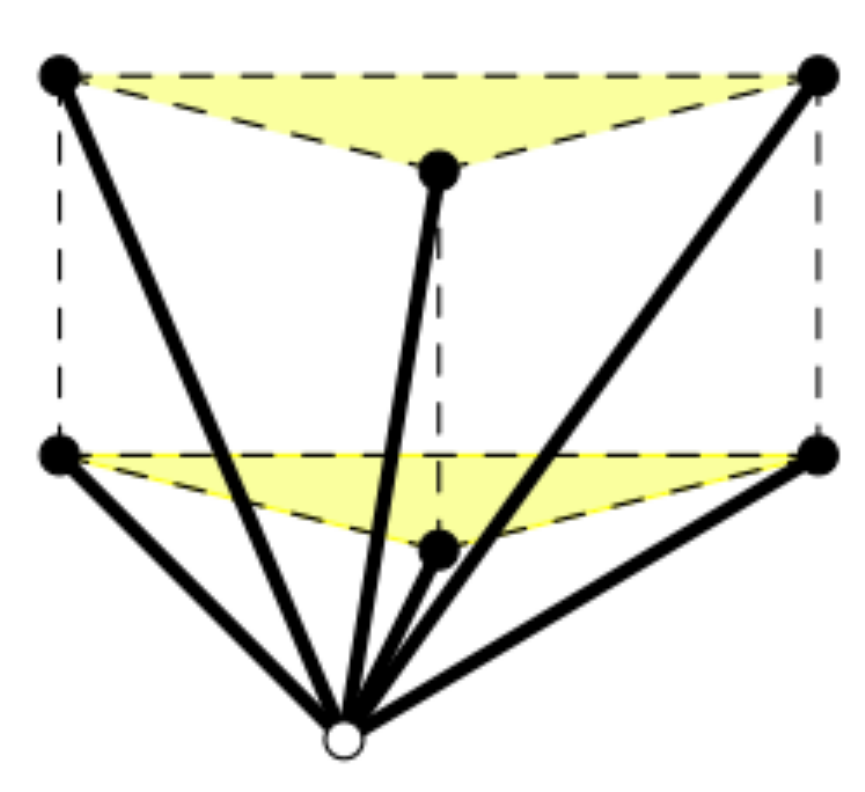}}
\vspace{-12pt}
 \caption{ Examples 1 and 2.}
 \label{FigDesigns12}
\vspace{-14pt}
\end{wrapfigure}

\medskip \noindent
{\bf Example 1.}\ The lattice of periods is generated by the edge-vectors of a tetrahedron. Thus, the vertices of the tetrahedron belong to the same orbit of framework vertices (black), but periods are not framework edges. We have $m=4$ edge orbits, with the four edge
representatives connecting a white vertex to the four black vertices of the tetrahedron, as shown in Fig.~\ref{FigDesigns12}(left).

\medskip \noindent
{\bf Example 2.}\ The lattice of periods is generated by the edge-vectors of a triangular prism. We have $m=6$ edge orbits, with the six edge representatives connecting a white vertex to the six black vertices of the triangular prism, as shown in Fig.~\ref{FigDesigns12}(right).

\begin{wrapfigure}{r}{0.55\textwidth}
\vspace{-20pt}
\centering
{\includegraphics[width=0.28\textwidth]{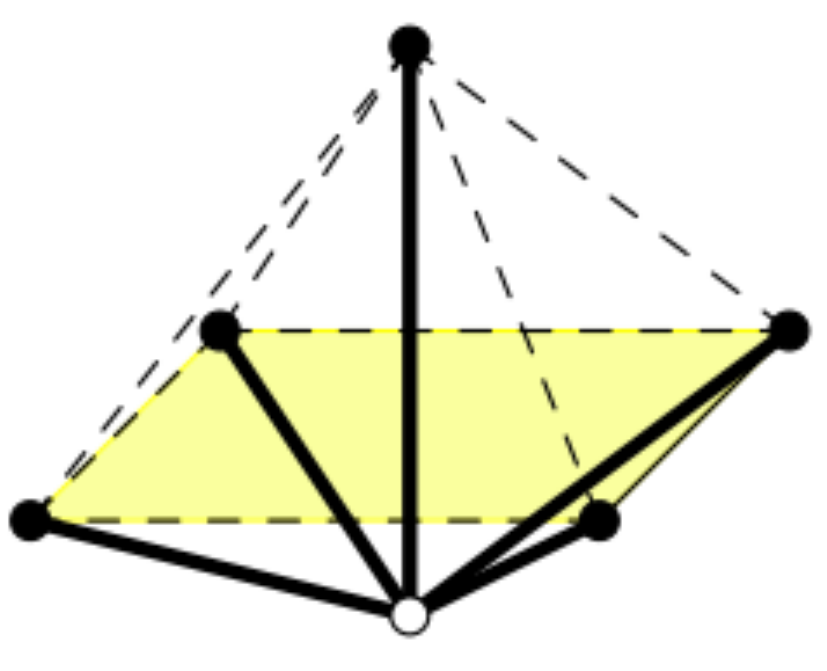}}
{\includegraphics[width=0.24\textwidth]{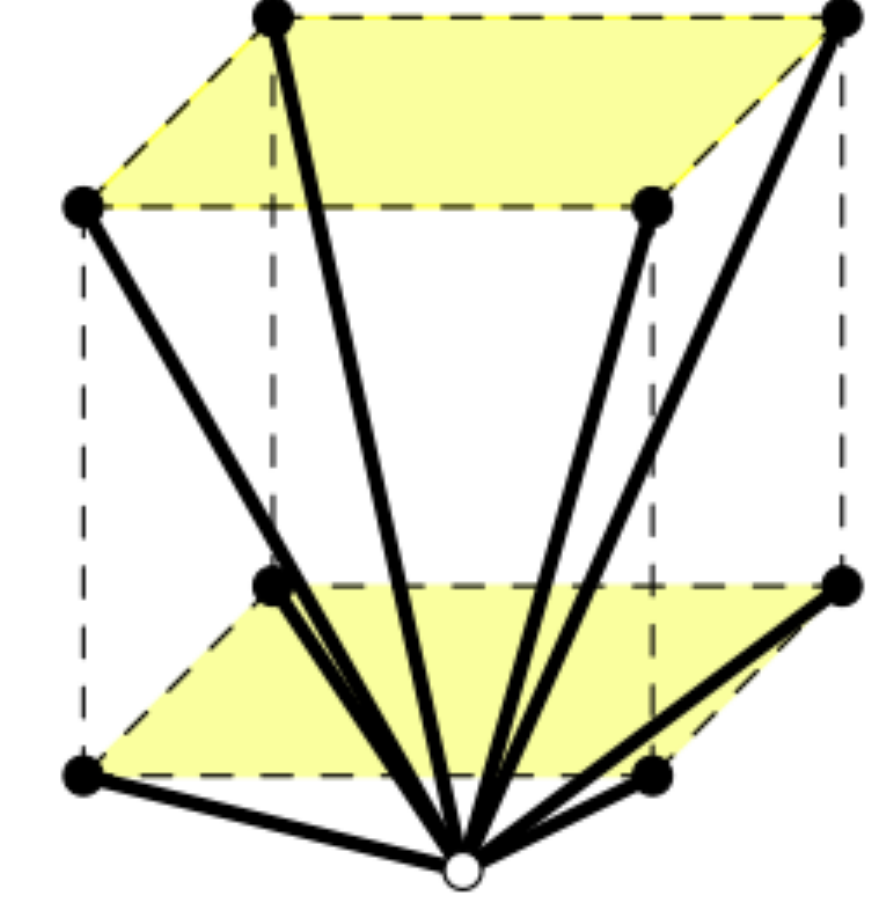}}
\vspace{-6pt}
 \caption{ Examples 3 and 4.}
 \label{FigDesigns34}
\vspace{-16pt}
\end{wrapfigure}

\medskip \noindent
{\bf Example 3.}\ The lattice of periods is generated by the edge-vectors of a regular pyramid with a squared basis. We have $m=5$ edge orbits, with the five edge representatives connecting a white vertex to the five black vertices of the pyramid, as shown in Fig.~\ref{FigDesigns34}(left). The reader will recognize this example as the one treated above in Section~\ref{sec:scenario}.

\medskip
\noindent
{\bf Example 4.}\ The lattice of periods is generated by the edge-vectors of a cube. We have $m=8$ edge orbits, with the eight edge representatives connecting a white vertex to the eight black vertices of the cube, as shown in Fig.~\ref{FigDesigns34}(right). This structure is discussed in more detail in \cite{BS7}. Auxetic capabilities were noticed earlier in \cite{Ma}, as cited in \cite{EL}.

\begin{wrapfigure}{l}{0.61\textwidth}
\vspace{-10pt}
\centering
{\includegraphics[width=0.32\textwidth]{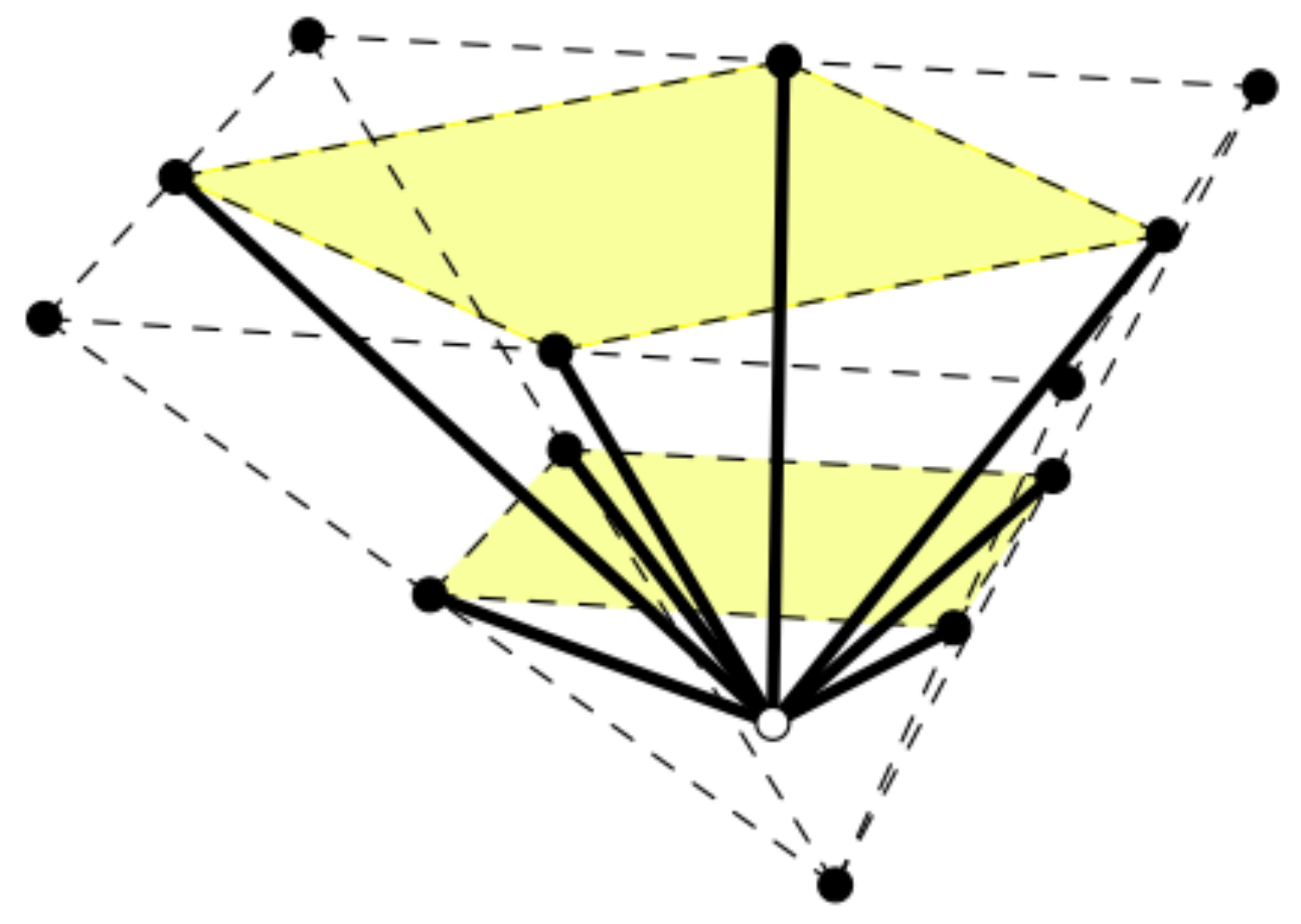}}
{\includegraphics[width=0.28\textwidth]{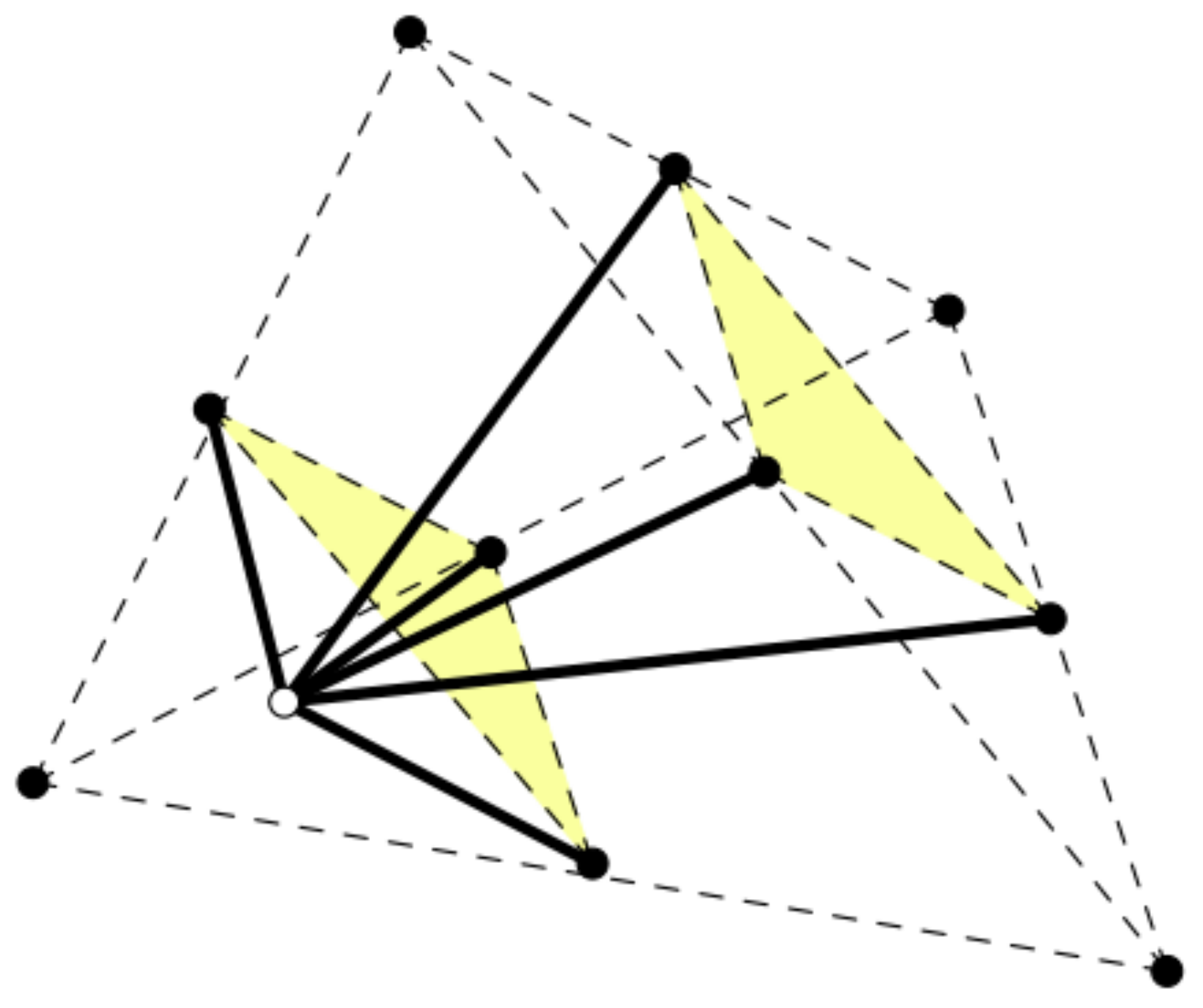}}
\vspace{-10pt}
 \caption{Examples 5 and 6.}
 \label{FigDesigns56}
\vspace{-16pt}
\end{wrapfigure}

\medskip \noindent
{\bf Example 5.}\ The lattice of periods is generated by the edge-vectors of a pyramid over a square. The square is seen as a `first ceiling'. We look at the doubled pyramid and call the doubled square the `second ceiling'. We have $m=8$ edge orbits, with the eight edge
representatives connecting a white vertex inside the pyramid to the four black vertices of the first ceiling and the four black midpoints of the boundary of the second ceiling, as shown in Fig.~\ref{FigDesigns56}(left).

\medskip \noindent
{\bf Example 6.}\ The lattice of periods is generated by the edge-vectors of a tetrahedron. We consider one face as a `first ceiling' and look at the doubled tetrahedron, where the doubled first ceiling provides the `second ceiling'. We have $m=6$ edge orbits, with the six edge
representatives connecting a white vertex inside the tetrahedron to the three black vertices of the first ceiling and to the black midpoints of the second ceiling, as shown in Fig.~\ref{FigDesigns56}(right). The corresponding structure in arbitrary dimension $d$ is discussed in \cite{BS7}.

\medskip \noindent 
Confirmation of auxetic capabilities with the methods of this paper is straightforward and is left to the reader.

\section{Conclusions}
\label{sec:conclusions} 

We have introduced a geometric theory of auxetic one-parameter deformations of periodic bar-and-joint frameworks, applicable in arbitrary dimension. Auxetic trajectories are characterized by the fact that the Gram matrix of a basis of periods evolves by keeping all tangent directions in the positive semidefinite cone. This is analogous to causal trajectories in special relativity, which have all tangents in the light cone. Thus, based on the geometry of the positive semidefinite cone, the infinitesimal auxetic cones can be determined and various auxetic vector fields or direction fields can be defined on the deformation space of any given periodic framework. Integral curves will be auxetic trajectories. 

\medskip
\noindent
For applications, dimensions two and three are most relevant. In dimension two, the structure of expansive periodic mechanisms is completely understood in terms of periodic pseudo-triangulations and an infinite range of auxetic designs follows directly from the stronger property of expansiveness. In dimension three or higher, expansive behavior is not yet sufficiently elucidated, but remains suggestive for auxetic design. However, auxetic capabilities may exist in the absence of expansive capabilities. 

\medskip 
\noindent
Higher dimensional considerations may prove important for exploring auxetic capabilities in
{\em quasi-crystals}, based on their description as projections of higher dimensional periodic
structures.

\medskip
\noindent
Last but not least, we observe that, by obviating the need for any actual determination of Poisson's ratios, our strictly geometric approach offers a precise and  rigorous method of auxetic investigations.


\end{document}